\documentclass[a4paper]{amsart}
\usepackage{amssymb}
\usepackage{amscd}
\usepackage{amsthm}
\usepackage{amsmath}
\usepackage{mathtools}
\usepackage{latexsym}
\usepackage{hyperref}
\usepackage[all]{xy}
\usepackage[utf8]{inputenc}
\usepackage{enumitem}

\usepackage{tikz-cd}

\usepackage{amsfonts}

\usepackage{booktabs}
\usepackage{multirow}
\usepackage{colortbl}

\usepackage{hyperref}
\usepackage{extarrows}
\usepackage{graphicx}
\usepackage{tikz}
\usepackage[table]{}

\theoremstyle{plain}
\newtheorem{theorem}{Theorem}
\newtheorem{corollary}[theorem]{Corollary}
\newtheorem{lemma}[theorem]{Lemma}
\newtheorem{proposition}[theorem]{Proposition}

\theoremstyle{definition}
\newtheorem{definition}[theorem]{Definition}
\newtheorem{setting}[theorem]{Setting}
\newtheorem{example}[theorem]{Example}

\theoremstyle{remark}

\newtheorem{remark}{Remark}

\numberwithin{theorem}{section}
\usepackage{newtxtext}
\usepackage{newtxmath}

\DeclareMathAlphabet\urwscr{U}{urwchancal}{m}{n}%
\DeclareMathAlphabet\rsfscr{U}{rsfso}{m}{n}
\DeclareMathAlphabet\euscr{U}{eus}{m}{n}
\DeclareFontEncoding{LS2}{}{}
\DeclareFontSubstitution{LS2}{stix}{m}{n}
\DeclareMathAlphabet\stixcal{LS2}{stixcal}{m} {n}

\newcommand{\rad}[1]{\operatorname{rad}(#1)}
\newcommand{\card}{\mbox{\rm{card\,}}}

\newcommand{\End}[2]{\mbox{\rm{End}}_{#1}(#2)}

\newcommand{\Soc}[1]{\mbox{\rm{Soc}}(#1)}

\newcommand{\supp}[1]{\mbox{\rm{supp}}(#1)}

\newcommand{\im}{\mbox{\rm{Im\,}}}
\newcommand{\Ker}{\mbox{\rm{Ker\,}}}

\newcommand{\Hom}[3]{\operatorname{Hom}_{#1}(#2,#3)}

\newcommand{\simp}[1]{\mbox{\rm{simp}--}{#1}}

\newcommand{\cs}[1]{\mbox{\rm{cs}}(#1)}

\newcommand{\Q}{\mathbb{Q}}   
\newcommand{\R}{\mathbb{R}}

\newcommand{\F}{\mathbb{F}}

\begin{document}

\title[Multiplicative bases and commutative semiartinian regular algebras]{Multiplicative bases and commutative semiartinian \\ von Neumann regular algebras}

\author{Kate\v{r}ina Fukov\'{a}}
\address{Charles University, Faculty of Mathematics
and Physics, Department of Algebra \\
Sokolovsk\'{a} 83, 186 75 Prague 8, Czech Republic}
\email{fukova@karlin.mff.cuni.cz}
\author{Jan Trlifaj}
\address{Charles University, Faculty of Mathematics
and Physics, Department of Algebra \\
Sokolovsk\'{a} 83, 186 75 Prague 8, Czech Republic}
\email{trlifaj@karlin.mff.cuni.cz}

\begin{abstract} Let $R$ be a semiartinian (von Neumann) regular ring with primitive factors artinian. The dimension sequence $\mathcal D _R$ is an invariant that captures the various skew-fields and dimensions occurring in the layers of the socle sequence of $R$. Though $\mathcal D _R$ does not determine $R$ up to an isomorphism even for rings of Loewy length $2$, we prove that it does so when $R$ is a commutative semiartinian regular $K$-algebra of countable type over a field~$K$. The proof is constructive: given the sequence $\mathcal D$, we construct the unique $K$-algebra of countable type $R = B_{\alpha,n}$ such that $\mathcal D = \mathcal D _R$ by a transfinite iterative construction from the base case of the $K$-algebra $R(\aleph_0,K)$ consisting of all eventually constant sequences in $K^{\aleph_0}$. Moreover, we prove that the $K$-algebras $B_{\alpha,n}$ possess conormed strong multiplicative bases despite the fact that the ambient $K$-algebras $K^{\kappa}$ do not even have any bounded bases for any infinite cardinal ~$\kappa$. 

Recently, a study of the number of limit models in AECs of modules \cite{M} has raised interest in the question of existence of strictly $\lambda$-injective modules for arbitrary infinite cardinals $\lambda$. In the final section, we construct examples of such modules over the $K$-algebra $R(\kappa,K)$ for each cardinal $\kappa \geq \lambda$.  
\end{abstract}

\date{\today}

\thanks{Research supported by GA\v CR 23-05148S, GAUK 101524 and SVV-2020-260721.}
    

\maketitle

\section{Introduction}

Let $R$ be a (right) semiartinian ring, that is, a ring $R$ such that the regular (right $R$-) module possesses a socle sequence $\mathcal S = ( S_\alpha \mid \alpha \leq \sigma + 1)$ with $S_{\sigma + 1} = R$. Then the ring $L_\sigma = R/S_\sigma$ is completely reducible, whence $L_\sigma$ is a finite ring direct product of full matrix rings over skew-fields. The ranks of these matrices and the skew-fields involved are invariants of $R$.

Assume moreover that $R$ is von Neumann regular with primitive factors artinian, cf. ~\cite[Chap.\ 6]{G}. Then similar invariants arise at each layer $L_\alpha$ ($\alpha \leq \sigma$) of $\mathcal S$. More precisely, for each $\alpha < \sigma$, $L_\alpha = S_{\alpha + 1}/S_\alpha$ is an infinite direct sum of rings without unit, each of which is a full matrix ring over a skew-field. Again, the ranks of these matrices and the skew-fields are invariants of $R$. Collected together over all $\alpha \leq \sigma$, these invariants form the \emph{dimension sequence} $\mathcal D _R$ of $R$ defined in ~\cite[Definition 2.2]{RTZ}, see also Lemma ~\ref{rtz-t} below.  

Isomorphic rings have equivalent dimension sequences (that is, they are factor equivalent in the sense of Definition ~\ref{facteq}), but the converse fails in general already for $\sigma = 1$ even for $R$ commutative and with countable layers (cf.\ Example ~\ref{non-unique}). However, if we restrict our setting to commutative regular semiartinian $K$-algebras $R$ that are of countable type, then the dimension sequence $\mathcal D _R$ determines $R$ up to a $K$-algebra isomorphism. This is the main result of our paper, Theorem ~\ref{unique'}. Here, $R$ is said to be of \emph{countable type} if its Loewy length $\sigma + 1$ is countable and for each $\alpha \leq \sigma$, the layer $L_\alpha$ is isomorphic, as a $K$-algebra without unit, to a direct sum of countably many copies of $K$.

Our proof of Theorem ~\ref{unique'} is constructive in the sense that for a given dimension sequence $\mathcal D$, the $K$-algebra $R$ such that $\mathcal D = \mathcal D_R$ is obtained by a transfinite iterative construction of rings based on ~\cite[Lemma 2.4]{EGT} which builds von Neumann regular semiartinian rings of larger Loewy length from sequences of rings of smaller lengths. The base case of the construction is the $K$-algebra $R(\aleph_0,K)$ consisting of all eventually constant sequences from ~$K^{\aleph_0}$.

\medskip
If $Q$ is a quiver without oriented cycles, then the standard $K$-basis $B$ of the finite dimensional path algebra $R = KQ$ is \emph{multiplicative} -- that is, $b b^\prime \in B \cup \{ 0 \}$ for all $b,b^\prime \in B$ -- and \emph{normed}, that is, $B$ contains a basis of each power of the Jacobson radical of $R$ as well as a complete set of pairwise orthogonal primitive idempotents of $R$. One of the highlights of classic representation theory of algebras, proved in \cite{BGRS}, says that also each finite dimensional $K$-algebra of finite representation type over an algebraically closed field $K$ possesses a normed multiplicative basis. 

Our Corollary ~\ref{multbs} shows that the $K$-algebras of countable type have a dual property: they possess conormed strong multiplicative bases. That is, if $R$ is of countable type, then $R$ has a multiplicative $K$-basis $B$ such that $b b^\prime \in B$ for all $b,b^\prime \in B$, and $B$ contains a basis of each term of the socle sequence of $R$. This contrasts with the fact, proved in Theorem \ref{fproducts}, that the ambient $K$-algebras $K^\kappa$ do not even possess any bounded bases for any infinite cardinal~$\kappa$. 

\medskip
As a by-product of our study, we prove in Theorem ~\ref{positive} that for each pair of infinite cardinals $\lambda \leq \kappa$, there exists a module $M_\lambda$ over the $K$-algebra $R = R(\kappa,K)$ such that $M_\lambda$ is strictly $\lambda$-injective (that is, $M_\lambda$ is $\lambda$-injective, but not $\lambda^+$-injective). This answers in the positive an old question of Paul Eklof, cf.\ \cite[Remark 2(iii)]{E}, reiterated recently in the context of AECs of modules by Mazari-Armida in ~\cite[Problem 3.9]{M}.

\section{Preliminaries}     

In what follows, $R$ we will denote an associative unital ring. Let $M$ be a (right $R$-) module. The \emph{socle sequence} of $M$ is defined as the strictly increasing sequence $\mathcal S = ( S_\alpha \mid \alpha \leq \tau )$ of submodules of $M$ satisfying $S_0 =0$, $S_{\alpha+1}/S_\alpha = \Soc{M/S_\alpha}$ for each $\alpha < \tau$, $S_\alpha = \bigcup_{\beta< \alpha} S_\beta$ for each limit ordinal $\alpha \leq \tau$, and $\Soc{M/S_\tau} = 0$. The ordinal $\tau$ as well as the sequence $\mathcal S$ are uniquely determined by $M$. If $S_\tau = M$ then $M$ is called \emph{semiartinian} and $\tau$ is the \emph{Loewy length} of $M$.

Notice that if $0 \neq M$ is finitely generated (e.g., when $M = R$) and $M$ is semiartinian, then the Loewy length of $M$ is a non-limit ordinal: $\tau = \sigma+1$ for some ordinal $\sigma$.

A ring $R$ is \emph{right semiartinian} if it is semiartinian as a right module over itself. Semiartinian rings are also known as Loewy rings and can equivalently be defined as rings over which each non-zero module $M$ satisfies $\Soc{M} \neq 0$. The socle sequence of the regular module $R_R$ is then called  \emph{the right socle sequence} of $R$.

A ring $R$ is said to have \emph{right primitive factors artinian} if $R/P$ is right artinian for each right primitive ideal $P$ of $R$; here, the term right \emph{primitive} refers to $P$ being the annihilator of a simple module.

A ring $R$ is called (von Neumann) \emph{regular} if for each $r \in R$ there exists $s_r \in R$ such that $rs_rr = r$. Moreover, $R$ is \emph{unit regular} if for each $r \in R$, $s_r$ can be taken to be a unit of $R$. An ideal $I$ of $R$ is regular if for each $i \in I$ there exists $j \in I$ such that $iji = i$.

Let $R$ be a right semiartinian regular ring with right primitive factors artinian, with a right socle sequence $\mathcal S$. As $R$ is regular, it is easy to see that $R$ is also left semiartinian, $\mathcal S$ is also its left socle sequence, and moreover, all left primitive factors of $R$ are artinian, cf. ~\cite[Theorem 6.2]{G}. Given this left-right symmetry, we will simply call such rings $R$ \emph{regular semiartinian rings with primitive factors artinian}. Notice that for each ideal $I \subsetneq R$, the ring $R/I$ is then also a regular semiartinian ring with primitive factors artinian.

\section{Dimension sequences for semiartinian regular rings}\label{dimseq}

We start by recalling a concrete description of regular semiartinian rings with primitive factors artinian from ~\cite[Theorem 2.1, Proposition 2.6]{RTZ}.

\begin{lemma}\label{rtz-t}
Let $R$ be a regular semiartinian ring with primitive factors artinian and $\mathcal S = (S_\alpha \mid \alpha \leq \sigma + 1)$ be the socle sequence of $R$.  

Then for each $\alpha \leq \sigma$ there are a cardinal $\lambda_\alpha > 0$, integers $n_{\alpha\beta} > 0$ ($\beta < \lambda_\alpha$) and skew-fields $K_{\alpha\beta}$ ($\beta < \lambda_\alpha$) such that the $\alpha$th layer $L_\alpha = S_{\alpha + 1}/S_\alpha$ of $R$ satisfies $L_\alpha \overset{\varphi_\alpha}\cong \bigoplus_{\beta < \lambda_\alpha} M_{n_{\alpha\beta}}(K_{\alpha\beta})$ as a ring without unit.

The cardinal $\lambda_\alpha$ is the number of homogenous components of $L_\alpha$, whence $\lambda_\alpha$ is infinite for $\alpha < \sigma$, and $\lambda_\sigma$ is finite.

The pre-image in $\varphi_\alpha$ of $M_{n_{\alpha\beta}}(K_{\alpha\beta})$ coincides with the $\beta$-th homogenous component $H_{\alpha \beta}$ of $L_\alpha$. This component is finitely generated for each $\beta < \lambda_\alpha$. Let $P_{\alpha \beta}$ denote a simple projective $R/S_\alpha$-module such that $H_{\alpha \beta} \cong P_{\alpha \beta}^{n_{\alpha\beta}}$. Then $K_{\alpha \beta} \cong \End{R/S_\alpha}{P_{\alpha \beta}} = \End{R}{P_{\alpha \beta}}$, and $\{ P_{\alpha \beta} \mid \alpha \leq \sigma, \beta < \lambda_\alpha \}$ is a representative set of all simple modules.
\end{lemma}

The structure of the rings $R$ from Lemma ~\ref{rtz-t} can be depicted as follows:

\begin{table}[!ht]
\begin{tabular}{ll|llll|l}
\cline{3-6}
\cellcolor[HTML]{34CDF9}{$\mathbf{L_\sigma}$}  &  & \cellcolor[HTML]{34CDF9}{$\mathbf{M_{n_{\sigma 0}}(K_{\sigma 0}) \, \oplus}$} & \cellcolor[HTML]{34CDF9}{\dots} & \cellcolor[HTML]{34CDF9}{$\mathbf{\oplus \, M_{n_{\sigma , \lambda_\sigma -1}}(K_{\sigma, \lambda_\sigma -1})}$} & \multicolumn{1}{l|}{\cellcolor[HTML]{34CDF9}\textbf{}} &                                                                                    \\ \cline{3-6}
\cellcolor[HTML]{96FFFB}{\dots} &  & \cellcolor[HTML]{96FFFB}{\dots} & \cellcolor[HTML]{96FFFB}{\dots} & \cellcolor[HTML]{96FFFB}{\dots} & \multicolumn{1}{l|}{\cellcolor[HTML]{96FFFB}\textbf{}} &                           \\ \cline{3-6}
\cellcolor[HTML]{38FFF8}{$\mathbf{L_\alpha}$}  &  & \cellcolor[HTML]{38FFF8}{$\mathbf{M_{n_{\alpha 0}}(K_{\alpha 0}) \, \oplus}$} & \cellcolor[HTML]{38FFF8}{\dots}  & \cellcolor[HTML]{38FFF8}{$\mathbf{\oplus \, M_{n_{\alpha\beta}}(K_{\alpha \beta}) \, \oplus}$ \dots} & \multicolumn{1}{l|}{\cellcolor[HTML]{38FFF8}{}} & $\beta < \lambda_\alpha$          \\ \cline{3-6}
\cellcolor[HTML]{68CBD0}{\dots} &  & \cellcolor[HTML]{68CBD0}{\dots} & \cellcolor[HTML]{68CBD0}{\dots} & \cellcolor[HTML]{68CBD0}{\dots} & \multicolumn{1}{l|}{\cellcolor[HTML]{68CBD0}\textbf{}} &                           \\ \cline{3-6}
\cellcolor[HTML]{00D2CB}{$\mathbf{L_1}$}  &  & \cellcolor[HTML]{00D2CB}{$\mathbf{M_{n_{10}}(K_{10}) \, \oplus}$} & \cellcolor[HTML]{00D2CB}{\dots} & \cellcolor[HTML]{00D2CB}{$\mathbf{\oplus \, M_{n_{1\beta}}(K_{1 \beta}) \, \oplus}$ \dots} & \multicolumn{1}{l|}{\cellcolor[HTML]{00D2CB}\textbf{}} & $\beta < \lambda_1$          \\ \cline{3-6}
\cellcolor[HTML]{329A9D}{$\mathbf{L_0 = Soc(R)}$}  &  & \cellcolor[HTML]{329A9D}{$\mathbf{M_{n_{00}}(K_{00}) \, \oplus}$}  & \cellcolor[HTML]{329A9D}{\dots} & \cellcolor[HTML]{329A9D}{$\mathbf{\oplus \, M_{n_{0 \beta}}(K_{0 \beta}) \, \oplus}$ \dots} & \multicolumn{1}{l|}{\cellcolor[HTML]{329A9D}\textbf{}} & $\beta < \lambda_0$ \\ \cline{3-6}
\end{tabular}
\end{table}

By Lemma ~\ref{rtz-t}, the sequence $\mathcal D _R = \{ ( \lambda_\alpha, \{ ( n_{\alpha\beta} , K_{\alpha\beta} ) \mid \beta < \lambda_\alpha \} ) \mid \alpha \leq \sigma \}$ is an invariant of $R$. It is called the \emph{dimension sequence} of $R$. The finite cardinal $\lambda_\sigma$ is called the \emph{top layer dimension} of $R$.

\begin{definition}\label{facteq} Let $R$ and $R^\prime$ be regular semiartinian rings with primitive factors artinian with the dimension sequences $\mathcal D _R = \{ ( \lambda_\alpha, \{ ( n_{\alpha\beta} , K_{\alpha\beta} ) \mid \beta < \lambda_\alpha \} ) \mid \alpha \leq \sigma \}$ and $\mathcal D_{R^\prime} = \{ ( \lambda^\prime_\alpha, \{ ( n^\prime_{\alpha\beta} , K^\prime_{\alpha\beta} ) \mid \beta < \lambda^\prime_\alpha \} ) \mid \alpha \leq \sigma^\prime \}$.

Then $R$ and $R^\prime$ are \emph{factor equivalent}, if $\sigma = \sigma^\prime$, $\lambda_\alpha = \lambda^\prime_\alpha$ for each $\alpha \leq \sigma$, and for each $\alpha \leq \sigma$ there exists a bijection $f_\alpha : \lambda_\alpha \rightarrow \lambda_\alpha$ such that $n_{\alpha\beta} = n^\prime_{\alpha f_\alpha (\beta)}$ and $K_{\alpha\beta} \cong K^\prime_{\alpha f_\alpha (\beta)}$, for all $\beta < \lambda_\alpha$.
\end{definition}

Clearly, isomorphic rings are factor equivalent, but factor equivalent rings need not be isomorphic even if they are commutative and of Loewy length $2$ (see Example ~\ref{non-unique} below).

\medskip
The following result from ~\cite{G} substantially simplifies our setting for commutative rings:

\begin{lemma}\label{compfa} (\cite[Theorem 3.2 and Corollary 4.2]{G}) Let $R$ be a commutative regular ring. Then $R$ is unit regular, and all primitive factors of $R$ are artinian.
\end{lemma}

Regular rings occur frequently in the setting of commutative semiartinian rings by an old result from ~\cite[Th\'eor\`eme 3.1]{NP}:

\begin{lemma}\label{nastpop} Let $R$ be any commutative semiartinian ring. Let $\rad{R}$ denote the Jacobson radical of $R$. Then the ring $R/\rad{R}$ is commutative, semiartinian, and regular.
\end{lemma}

Let us now characterize the commutativity of the factor-rings $R/S_\alpha$ for each $\alpha \leq \sigma$.

\begin{lemma}\label{relcommutat} Let $R$ be semiartinian regular ring with primitive factors artinian and let $\mathcal D_R = \{ ( \lambda_\alpha, \{ ( n_{\alpha\beta} , K_{\alpha\beta} ) \mid \beta < \lambda_\alpha \} ) \mid \alpha \leq \sigma \}$.

Then for each $\alpha \leq \sigma$, the following conditions are equivalent:
\begin{enumerate}
\item[(i)] $R/S_\alpha$ is commutative,
\item[(ii)] For each $\beta < \lambda_\alpha$, $K_{\alpha\beta}$ is commutative and $n_{\alpha\beta} = 1$,
\item[(iii)] For each $\beta < \lambda_\alpha$, $K_{\alpha\beta}$ is commutative, and $R/S_\alpha$ is isomorphic to a subring $\bar{R}_\alpha$ of the ring direct product $\prod_{\beta < \lambda_\alpha} K_{\alpha \beta}$ such that $\Soc{\bar{R}_\alpha} = \bigoplus_{\beta < \lambda_\alpha} K_{\alpha \beta}$.
\end{enumerate}
\end{lemma}
\begin{proof} (i) implies (ii). Since $R/S_\alpha$ is commutative, $n_{\alpha\beta} = 1$ for all $\beta < \lambda_\alpha$. Let $M_{\alpha \beta}/S_{\alpha}$ be the maximal ideal of the ring $R/S_\alpha$ such that $R/S_\alpha = M_{\alpha \beta}/S_{\alpha} \oplus P_{\alpha \beta}$. Then $R/M_{\alpha \beta}$ is a field isomorphic to $\End{R}{P_{\alpha \beta}}$, whence the skew-field $K_{\alpha \beta} \cong \End{R/S_\alpha}{P_{\alpha \beta}} = \End{R}{P_{\alpha \beta}}$ is commutative.

(ii) implies (iii). We have $\Soc{R/S_\alpha} = \bigoplus_{\beta < \lambda_\alpha} P_{\alpha \beta}$. Since $n_{\alpha \beta} = 1$, for each $\beta < \lambda_\alpha$, $P_{\alpha \beta}$ is the $\beta$th homogenous component of $\Soc{R/S_\alpha}$, and it is a simple direct summand in $R/S_\alpha$. Since $P_{\alpha \beta}$ is projective as an $R/S_\alpha$-module, there is a maximal right ideal $M_{\alpha \beta}$ of $R/S_\alpha$ such that $R/S_\alpha = M_{\alpha \beta} \oplus P_{\alpha \beta}$. Moreover, $M_{\alpha \beta} \supseteq \bigoplus_{\beta \neq \gamma < \lambda_\alpha} P_{\alpha \gamma}$, and $M_{\alpha \beta}$ contains no simple submodule isomorphic to $P_{\alpha \beta}$. Being a homogenous component of $\Soc{R/S_\alpha}$, $P_{\alpha \beta}$ is a two-sided ideal in $R/S_\alpha$, and since $P_{\alpha \beta}$ is projective, also $M_{\alpha \beta}$ is a two-sided ideal in $R/S_\alpha$. That is, $R = M_{\alpha \beta} \boxplus P_{\alpha \beta}$ is a ring direct sum decomposition of $R/S_\alpha$. Let $e_{\alpha \beta}$ be the central idempotent in $R/S_\alpha$ such that $P_{\alpha \beta} = e_{\alpha \beta} (R/S_\alpha)$.    

Let $M = \bigcap_{\beta < \lambda_\alpha} M_{\alpha \beta}$. Then $M \cap \Soc{R/S_\alpha} = 0$, whence $M = 0$. Thus, the ring homomorphism $\varphi : R/S_\alpha \to \prod_{\beta < \lambda_\alpha} R/M_{\alpha \beta}$ defined by $\varphi(r) = (r + M_{\alpha \beta})_{\beta < \lambda_\alpha}$ is monic. Moreover, $\varphi (\Soc{R/S_\alpha}) = \varphi (\bigoplus_{\beta < \lambda_\alpha} P_{\alpha \beta}) = \bigoplus_{\beta < \lambda_\alpha} R/M_{\alpha \beta}$. The ring isomorphisms $K_{\alpha \beta} \cong \End{R/S_\alpha}{P_{\alpha \beta}} \cong e_{\alpha \beta} (R/S_\alpha) e_{\alpha \beta} = e_{\alpha \beta} (R/S_\alpha) \cong R/M_{\alpha \beta}$ ($\beta < \lambda_\alpha$) then yield the commutativity of $K_{\alpha \beta}$ as well as the final claim for $\bar{R}_\alpha = \varphi(R/S_\alpha)$.  

(iii) implies (i). This is clear since in the setting of (iii), the ring $\prod_{\beta < \lambda_\alpha} K_{\alpha \beta}$ is commutative.
\end{proof}

It follows that commutativity can be tested at the bottom level of the dimension sequence.

\begin{corollary}\label{commutat}
Let $R$ be a semiartinian regular ring with primitive factors artinian and let $\mathcal D_R = \{ ( \lambda_\alpha, \{ ( n_{\alpha\beta} , K_{\alpha\beta} ) \mid \beta < \lambda_\alpha \} ) \mid \alpha \leq \sigma \}$. Then the following conditions are equivalent:
\begin{enumerate}
\item[(i)] $R$ is commutative,
\item[(ii)] For all $\alpha \leq \sigma$ and $\beta < \lambda_\alpha$, $K_{\alpha\beta}$ is commutative and $n_{\alpha\beta} = 1$.
\item[(iii)] For each $\beta < \lambda_0$, $K_{0 \beta}$ is commutative and $n_{0 \beta} = 1$.
\end{enumerate}
\end{corollary}

Notice that part (iii) of Lemma ~\ref{relcommutat} restricts the size of $R$ and of the fields $K_{\alpha \beta}$, since $R$ is a subring of $\prod_{\beta < \lambda_0} K_{0 \beta}$. For further restrictions on possible parameters of dimension sequences, we refer to ~\cite{Z1} and ~\cite{Z2}.  

\medskip
Our main interest will be in the following particular instance of the commutative setting.

\begin{definition}\label{onefield} Let $K$ be a field. We will denote by $\mathfrak R _K$ the class of all $K$-algebras $R$, which are commutative semiartinian regular rings, such that for their dimension sequences $\mathcal D_R = \{ ( \lambda_\alpha, \{ ( n_{\alpha\beta} , K_{\alpha\beta} ) \mid \beta < \lambda_\alpha \} ) \mid \alpha \leq \sigma \}$ the following holds: $K_{\alpha \beta} = K$ for all $\alpha \leq \sigma$ and $\beta < \lambda_\alpha$, and the isomorphisms $\varphi_\alpha$ ($\alpha \leq \sigma$) from Lemma ~\ref{rtz-t} are $K$-linear.
\end{definition}
Equivalently, $\mathfrak R _K$ is the class of all commutative semiartinian regular $K$-algebras $R$ such that $K$ is the endomorphism ring of all the simple ($R$- or $R/S_\alpha$-) modules $P_{\alpha \beta}$ ($\alpha \leq \sigma$ and $\beta < \lambda_\alpha$), that is, all the endomorphisms are given by multiplications by the elements of $K$.

Thus, if $\simp R$ denotes a representative set of all simple modules, then for each $S \in \simp R$ there is a unique $\alpha \leq \sigma$ such that $S$ embeds into $R/S_\alpha$, the $S$-homogenous component $H_S$ of $S_{\alpha +1}/S_\alpha$ is simple, that is, $H_S \cong S$, and the endomorphisms of $S$ are given by multiplications by the elements of $K$.

A $K$-algebra $R \in \mathfrak R _K$ can be depicted as follows:

\vfill\eject

\medskip
\begin{table}[h]
\begin{tabular}{ll|llll|l}
\cline{3-6}
\cellcolor[HTML]{34CDF9}{$\mathbf{L_\sigma \overset{\varphi_\sigma}\cong K^{(\lambda_{\sigma})}}$}  &  & \cellcolor[HTML]{34CDF9}{$\mathbf{K \, \oplus}$} & \cellcolor[HTML]{34CDF9}{\dots} & \cellcolor[HTML]{34CDF9}{$\mathbf{\oplus \, K}$} & \multicolumn{1}{l|}{\cellcolor[HTML]{34CDF9}\textbf{}} &                                                                                  \\ \cline{3-6}
\cellcolor[HTML]{96FFFB}{\dots} &  & \cellcolor[HTML]{96FFFB}{\dots} & \cellcolor[HTML]{96FFFB}{\dots} & \cellcolor[HTML]{96FFFB}{\dots} & \multicolumn{1}{l|}{\cellcolor[HTML]{96FFFB}\textbf{}} &                           \\ \cline{3-6}
\cellcolor[HTML]{38FFF8}{$\mathbf{L_\alpha \overset{\varphi_\alpha}\cong K^{(\lambda_{\alpha})}}$}  &  & \cellcolor[HTML]{38FFF8}{$\mathbf{K \, \oplus}$} & \cellcolor[HTML]{38FFF8}{\dots}  & \cellcolor[HTML]{38FFF8}{$\mathbf{\oplus \, K \, \oplus}$ \dots} & \multicolumn{1}{l|}{\cellcolor[HTML]{38FFF8}{}} &           \\ \cline{3-6}
\cellcolor[HTML]{68CBD0}{\dots} &  & \cellcolor[HTML]{68CBD0}{\dots} & \cellcolor[HTML]{68CBD0}{\dots} & \cellcolor[HTML]{68CBD0}{\dots} & \multicolumn{1}{l|}{\cellcolor[HTML]{68CBD0}\textbf{}} &                           \\ \cline{3-6}
\cellcolor[HTML]{00D2CB}{$\mathbf{L_1 \overset{\varphi_1}\cong K^{(\lambda_{1})}}$}  &  & \cellcolor[HTML]{00D2CB}{$\mathbf{K \, \oplus}$} & \cellcolor[HTML]{00D2CB}{\dots} & \cellcolor[HTML]{00D2CB}{$\mathbf{\oplus \, K \, \oplus}$ \dots} & \multicolumn{1}{l|}{\cellcolor[HTML]{00D2CB}\textbf{}} &           \\ \cline{3-6}
\cellcolor[HTML]{329A9D}{$\mathbf{L_0 = Soc(R) \overset{\varphi_o}\cong K^{(\lambda_{0})}}$}  &  & \cellcolor[HTML]{329A9D}{$\mathbf{K \, \oplus}$}  & \cellcolor[HTML]{329A9D}{\dots} & \cellcolor[HTML]{329A9D}{$\mathbf{\oplus \, K \, \oplus}$ \dots} & \multicolumn{1}{l|}{\cellcolor[HTML]{329A9D}\textbf{}} &  \\ \cline{3-6}
\end{tabular}
\end{table}

\begin{definition}      
If $R \in \mathfrak R _K$, the ordinal $\sigma$ is countable, and so are the cardinals $\lambda_\alpha$ ($\alpha \leq \sigma$), then $R$ will be called of \emph{countable type}.
\end{definition}
Notice that in this case $\lambda_\alpha = \aleph_0$ for all $\alpha < \sigma$, and $K_{\alpha \beta} = K$ for all $\beta < \lambda_\alpha$ and $\alpha \leq \sigma$. So the ordinal $\sigma$ and the finite number $\lambda_\sigma$ are the only parameters in the dimension sequence of a $K$-algebra $R$ of countable type that can vary. Since all ideals of $K$-algebras of countable type are countably generated, these algebras are hereditary by ~\cite[Proposition 2.14]{G}.

It easily follows from the definitions above that if $K$ is a field and $R_1, R_2 \in \mathfrak R _K$, then also the direct sum $R = R_1 \boxplus R_2 \in \mathfrak R _K$. If moreover $R_1$ and $R_2$ are of countable type, then so is $R$.

The following lemma shows that when investigating the structure of rings in the class $\mathfrak R _K$, we can restrict ourselves to rings of the top layer dimension $1$:

\begin{lemma}\label{topone'} Let $R \in \mathfrak R _K$ be of Loewy length $\sigma + 1$ and of the top layer dimension $\lambda_\sigma = n$. Then there is a $K$-algebra direct product decomposition $R = R_0 \boxplus \dots \boxplus R_{n-1}$ where for each $i < n$, $R_i \in \mathfrak R _K$ has Loewy length $\sigma + 1$ and top layer dimension $1$.  
\end{lemma}
\begin{proof} By assumption, we have a $K$-algebra isomorphism $R/S_\sigma = L_\sigma \overset{\varphi_\sigma}\cong K^n$. Let $\{ e_i \mid i < n \}$ be the canonical basis of $K^n$. By ~\cite[Proposition 2.18]{G}, there exists a complete set of orthogonal idempotents $\{ f_i \mid i < n \}$ in $R$ such that $e_i = \varphi_{\sigma}(f_i + S_{\sigma})$ for each $i < n$. Since $R$ is commutative, we have the $K$-algebra finite direct product decomposition $R = f_0R \boxplus \dots \boxplus f_{n-1}R$, and hence also $S_\alpha = \bigoplus_{i < n} f_iS_{\alpha}$ for each $\alpha \leq \sigma$. For each $i < n$, $f_i \in R \setminus S_\sigma$ is an idempotent whence  for each $\alpha \leq \sigma$, $f_iS_\alpha = f_iR \cap S_\alpha$, and $f_iS_{\alpha+1}/f_iS_\alpha \cong (f_iS_{\alpha+1} + S_\alpha)/S_\alpha = f_i L_\alpha$.

Thus, for each $\alpha \leq \sigma$, there is an isomorphism of $K$-algebras without unit $$\bigoplus_{i < n} f_iS_{\alpha+1}/f_iS_\alpha \overset{\phi_\alpha}\cong L_\alpha,$$ and $( f_iS_{\alpha} \mid \alpha \leq \sigma + 1 )$ is the socle sequence of $f_iR$ for each $i < n$. It follows that for each $i < n$, $f_iR \in \mathfrak R _K$.

Finally, $f_iR/f_iS_{\sigma} \overset{\rho}\cong (f_iR + S_\sigma)/S_\sigma \overset{\pi}\cong e_i (R/S_\sigma) \cong K$, where $\rho$ is the restriction of $\phi_\sigma$ to $f_iR/f_iS_{\sigma}$ and $\pi$ the restriction of $\varphi_\sigma$ to $(f_iR + S_\sigma)/S_\sigma$. Hence, the top layer dimension of $f_iR$ is $1$.
\end{proof}

\begin{remark}\label{remdec'} Notice that if $R \in \mathfrak R _K$ is of Loewy length $\sigma + 1$ and of the top layer dimension $1$, then the $K$-algebra isomorphism $\varphi_\sigma : R/S_\sigma \to K$ yields a $K$-vector space decomposition $R = S_\sigma \oplus 1_R \cdot K$ where $1_R$ is the unit of $R$.
\end{remark}
 
\section{Ring semidirect products}\label{semidir}

Extensions and split extensions are among the key concepts of the structure theory of modules. We will employ analogous, but less well-known, concepts of an extension and a semidirect product of rings:  

\begin{definition}\label{semid}
Let $R$ and $S$ be rings, and $I$ a two-sided ideal in $R$.

\begin{enumerate}
\item[(i)] $R$ is an \emph{extension} of $I$ by $S$ provided that there exists a surjective ring homomorphism $\pi : R \to S$ such that $\Ker \pi = I$. That is, there is a short exact sequence $$\mathcal E : 0 \to I \to R \overset{\pi}\to S \to 0.$$
Two extensions $\mathcal E$ and $\mathcal E ^\prime : 0 \to I \to R^\prime \overset{\pi^\prime}\to S \to 0$ of $I$ by $S$ are said to be \emph{equivalent} in case there is a ring homomorphism $\phi : R \to R^\prime$ such that the following diagram is commutative
$$\begin{CD}
0@>>>  I@>>> R@>{\pi}>>  S@>>>  0\\
@.     @|    @V{\phi}VV @|    @.\\
0@>>>  I@>>> {R^\prime}@>{\pi^\prime}>>  S@>>>  0
\end{CD}$$
\item[(ii)] $R$ is a \emph{semidirect product} of $I$ and $S$ (or a \emph{split extension} of $I$ by $S$) in case there exists a split ring homomorphism $\pi : R \to S$ such that $\Ker \pi = I$. Here, $\pi$ is called \emph{split} in case there is a ring homomorphism $\nu : S \to R$ such that $\pi \nu = \hbox{id}_S$.
\end{enumerate}
\end{definition}

We will often use the following simple fact (cf. ~\cite[Lemma 1.3]{G}): If $R$ is an extension of $I$ by $S$ in the sense of Definition ~\ref{semid}(i), then $R$ is regular, iff both $I$ and $R/I$ are regular.

\medskip
There is an easy alternative description of semidirect products:

\begin{lemma}\label{alter} Let $R$ and $S$ be rings, and $I$ a two-sided ideal in $R$. Then $R$ is a semidirect product of $I$ and $S$, iff there exists an idempotent ring endomorphism $\epsilon$ of $R$ such that $\Ker \epsilon = I$ and $\im \epsilon \cong S$.
\end{lemma}
\begin{proof} Assume that $\pi$ is a ring homomorphism $\pi : R \to S$ such that $\Ker \pi = I$, and $\nu$ is a ring homomorphism $\nu : S \to R$ such that $\pi \nu = \hbox{id}_S$. Then $\epsilon = \nu \pi$ is an idempotent ring endomorphism of $R$ with $\Ker \epsilon = \Ker \pi = I$ and $\im \epsilon \cong \im \pi = S$.   

Conversely, let $\epsilon$ be an idempotent ring endomorphism of $R$ such that $\Ker \epsilon = I$ and there is a ring isomorphism $\psi: \im \epsilon \cong S$. Let $\mu : \im \epsilon \hookrightarrow R$ denote the inclusion, and let $\pi = \psi \epsilon : R \to S$. Then $\Ker \pi = \Ker \epsilon = I$, and for $\nu = \mu \psi ^{-1}$ we have $\pi \nu = \psi \epsilon \mu  \psi ^{-1} = \hbox{id}_S$ because $\epsilon$ is an idempotent.   
\end{proof}

Clearly, if $R$ is a semidirect product of $I$ and $S$, then $R$ is an extension of $I$ by $S$, but the converse fails in general even if $S$ is a finite field - just consider the extension/short exact sequence $0 \to {p\mathbb Z} \to \mathbb Z \to S = \mathbb Z/p{\mathbb Z} \to 0$ for a prime integer $p$.

We will be particularly interested in extensions fitting in the following particular setting:
    
\begin{setting}\label{set} Let $K$ be a field, $\kappa$ an infinite cardinal, and $\mathcal R = ( R_\alpha \mid \alpha < \kappa)$ a sequence of $K$-algebras. Let $P = \prod_{\alpha < \kappa} R_\alpha$ denote the $K$-algebra which is the product of the algebras in $\mathcal R$, and let $I = \bigoplus_{\alpha < \kappa} R_\alpha$ (so $I$ is an ideal in $P$). Let $R$ be an extension of $I$ by $K$, so $I$ is an ideal of $R$, and there is a ring isomorphism $\varphi : R/I \cong K$.

If moreover $R$ is a $K$-algebra and the ring isomorphism $\varphi$ is $K$-linear, then $R$ is isomorphic to the $K$-subalgebra $S = I \oplus 1_P \cdot K$ of $P$. We will denote this $K$-subalgebra by $R(\kappa,K,\mathcal R)$. Clearly, $R(\kappa,K,\mathcal R)$ is a semidirect product of $I = \bigoplus_{\alpha < \kappa} R_\alpha$ and $K$. 
\end{setting}

\begin{remark}\label{comments1} Notice that if $\mathcal R ^\prime = ( R^\prime _\alpha \mid \alpha < \kappa )$ is a sequence of $K$-algebras such that $R^\prime _\alpha \cong R_\alpha$ as $K$-algebras for all $\alpha < \kappa$, then also $R(\kappa,K,\mathcal R) \cong R(\kappa,K,\mathcal R ^\prime)$ as $K$-algebras. 

Hence, in the definition of $R(\kappa,K,\mathcal R)$ above, we can replace $\kappa$ by any infinite indexing set $A$. In this case, we will use the notation $R(A,K,\mathcal R)$ rather than $R(\kappa,K,\mathcal R)$.
\end{remark}

We will also employ the following uniqueness property of the algebra $R(\kappa,K,\mathcal R)$:

\begin{lemma}\label{extensions} Let $K$ be a field and $\kappa$ an infinite cardinal. Let $\mathcal R = ( R_\alpha \mid \alpha < \kappa )$ be a sequence of $K$-algebras, and $I = \bigoplus_{\alpha < \kappa} R_\alpha$ be the ideal of the $K$-algebra $P = \prod_{\alpha < \kappa} R_\alpha$ from Setting ~\ref{set}.

Let $R$ be a $K$-algebra containing an ideal $J$ which is, as a $K$-algebra without unit, isomorphic to $I$. Moreover, assume that there is a $K$-algebra isomorphism $R/J \cong K$. Then the $K$-algebras $R$ and $R(\kappa,K,\mathcal R)$ are isomorphic.
\end{lemma}
\begin{proof} By assumption, there is an isomorphism of $K$-algebras without unit $\phi : J \to I$. Since the $K$-algebra $R/J$ has dimension $1$, $R$ decomposes as $R = J \oplus 1_R \cdot K$. By the definition of $R(\kappa,K,\mathcal R)$, we have $R(\kappa,K,\mathcal R) = I \oplus 1_P \cdot K$. Define $\psi : R \to R(\kappa,K,\mathcal R)$ by the formula $\psi (j + 1_R\cdot k) = \phi(j) + 1_P \cdot k$ for all $j \in J$ and $k \in K$. Then $\psi$ is a $K$-algebra isomorphism of $R$ onto $R(\kappa,K,\mathcal R)$.
\end{proof}

\begin{corollary}\label{decomp} In the Setting \ref{set}, let $\beta < \kappa$ and $A = \kappa \setminus \{ \beta \}$. Then there is a $K$-algebra isomorphism $R(\kappa,K,\mathcal R) \cong R(A,K,\mathcal R) \boxplus R_{\beta}$.
\end{corollary}  
\begin{proof} This follows from Lemma ~\ref{extensions} for $R = R(A,K,\mathcal R) \boxplus R_{\beta}$ and $J = \bigoplus_{\alpha \in A} R_\alpha \oplus R_\beta$, since $J \cong I$ and $R = J \oplus 1_R \cdot K$.
\end{proof}

\medskip
Unlike the general case, if $K$ is a finite field, then any extension in Setting ~\ref{set} is necessarily split:  

\begin{proposition}\label{finite} Assume that $K$ is finite in the Setting ~\ref{set}. Then $R$ is a semidirect product of $I$ and $K$.  
\end{proposition}
\begin{proof} We have a short exact sequence $\mathcal E : 0 \to I \to R \overset{\pi}\to K \to 0$ where $\pi(r) = \varphi(r + I)$, so $\Ker \pi = I$. Let $\mu : K \to R$ be any lifting of $\pi$, that is, a map such that $\pi(\mu(k)) = k$ for each $k \in K$. W.l.o.g., $\mu(1_{K}) = 1_R (= 1_P)$ and $\mu(0_{K}) = 0_R (= 0_P)$. Let $x, y \in K$. Then $\pi(\mu(x+y) - \mu(x) - \mu(y)) = 0$ and $\pi(\mu(x\cdot y) - \mu(x)\mu(y)) = 0$ whence $\mu(x)+\mu(y) - \mu(x+y) \in I$ and $\mu(x)\cdot \mu(y) - \mu(x\cdot y) \in I$. In particular, there is a finite subset $A_{(x,y)} \subseteq \kappa$ such that for each $\alpha \in \kappa \setminus A_{(x,y)}$, the $\alpha$th component of $\mu(x+y)$ coincides with the sum of the $\alpha$th components of $\mu(x)$ and $\mu(y)$, and the $\alpha$th component of $\mu(x\cdot y)$ coincides with the product of the $\alpha$th components of $\mu(x)$ and $\mu(y)$. Let $A = \bigcup_{(x,y) \in K \times K} A_{(x,y)}$. Since $K$ is finite, $A$ is a finite subset of $\kappa$. Define $\nu : K \to R$ as follows: for each $k \in K$, the $\alpha$th component of $\nu(k)$ equals the $\alpha$th component of $\mu(k)$ for each $\alpha \in \kappa \setminus A$, and the $\alpha$th component of $\nu(k)$ equals $k$ for each $\alpha \in A$. Then $\nu$ is a ring homomorphism from $K$ to $R$ and $\im{(\mu - \nu)} \subseteq I$, whence $\pi\nu = \pi\mu = \hbox{id}_{K}$ proving that $R$ is a semidirect product of $I$ and $K$.
\end{proof}

The following simple example also fits Setting ~\ref{set}. It will serve as a base case for more complex constructions later on.

\begin{example}\label{exkappa} Let $K$ be a field, $\kappa$ an infinite cardinal, and $\mathcal R = ( R_\alpha \mid \alpha < \kappa)$ be such that $R_\alpha$ is just a field $K$ for each $\alpha < \kappa$. Then the $K$-algebra $R = R(\kappa,K,\mathcal R)$ is the $K$-subalgebra of $P = K^\kappa$ consisting of all sequences from $P$ that are constant except for finitely many terms. Since the ideal $I = K^{(\kappa)}$ of $R$ is regular, so is $R$. Notice that $R$ is commutative semiartinian of Loewy length $2$, $R \in \mathfrak R _K$ (and if moreover $\kappa = \aleph_0$, then $R$ is of countable type). We will also use the simplified notation $R = R(\kappa,K)$.
\end{example}
 
\medskip
For any two modules $A$ and $B$, all split extensions of $A$ by $B$ are equivalent, and hence their middle terms are isomorphic. This is not true of semidirect products of rings in general, even in our particular Setting ~\ref{set}:

\begin{example}\label{s-semid} Let $K$ be a field of characteristic $p \neq 0$. Let $\varphi$ denote the Frobenius endomorphism of $K$ (i.e., $\varphi(k) = k^p$ for all $k \in K$). Let $\kappa$ be an infinite cardinal. Further, we let $R_\alpha = K$ for each $\alpha < \kappa$. Then $I = \bigoplus_{\alpha < \kappa} R_\alpha = K^{(\kappa)}$ is an ideal in the $K$-algebra $P = \prod_{\alpha < \kappa} R_\alpha = K^\kappa$. For each $\alpha < \kappa$, we denote by $e_\alpha$ the primitive idempotent in $P$ such that $e_\alpha P = R_\alpha$.

Define $\nu : K \to P$ as follows: for each $k \in K$, $\nu(k)$ is the element of $P$ whose all components equal $k^p$. Since $K$ has characteristic $p$, $\nu$ is a ring monomorphism. Let $R = I \oplus \nu(K)$. Then $R$ is a subring of $P$ (but not a $K$-subalgebra of $P$ unless $k^p = k$ for all $k \in K$, i.e., unless $K = \mathbb Z _p$), and $R$ is regular. There is an extension $\mathcal E : 0 \to I \to R  \overset{\pi}\to K \to 0$ where the ring homomorphism $\pi$ is defined by $\pi(r) = k$ for each $r \in R$ of the form $r = i + \nu(k)$, where $i \in I$ and $k \in K$. So we are again in Setting ~\ref{set}. Moreover, since $\pi \nu = \hbox{id}_K$, $R$ is a semidirect product of $I$ and $K$.

Notice that $K$ is a $\nu(K)$-algebra, hence so is $P$, and $R$ is a $\nu(K)$-subalgebra of $P$. As the fields $K$ and $\nu(K)$ are isomorphic, $R$ is also a $K$-algebra.  Clearly, $R$ is a subring of $R(\kappa,K)$ and $E = \{ e_\alpha \mid \alpha < \kappa \}$ is the set of all primitive idempotents in both regular rings $R$ and $R(\kappa,K)$. In particular, $I = K^{(\kappa)} = \Soc R = \Soc {R(\kappa,K)}$.

Let define the map $\psi : R(\kappa,K) \to R$ by $\psi(\sum_{i < m} e_{\alpha_i} k_{\alpha_i} + 1_P\cdot k) = \sum_{i < m} e_{\alpha_i} \nu(k_{\alpha_i}) + 1_P \nu(k)$. Since $\nu$ is a ring monomorphism, so is $\psi$. We will distinguish two cases.

\medskip
(i) \emph{$K$ is a perfect field} (i.e. $\varphi$ is an automorphism of $K$). Then $\im \nu = 1_P\cdot K$, so $\psi$ is surjective, and the rings $R(\kappa,K)$ and $R$ are isomorphic.

However, the extensions $\mathcal E$ and $\mathcal E ^\prime : 0 \to I \to R(\kappa,K) \overset{\pi^\prime}\to K \to 0$ where $\pi^\prime (i + 1_p\cdot k) = k$ are not equivalent unless $\varphi = \hbox{id}_K$ (i.e., unless $K = \mathbb Z _p$). Indeed, any ring isomorphism $\phi$ of $R(\kappa,K)$ on to $R$ that is identity on $I$ must satisfy $e_\alpha k = \phi (e_\alpha k) = e_\alpha \phi (1_P\cdot k)$ for all $\alpha < \kappa$ and $k \in K$. However, if $k^p \neq k$ for some $k \in K$, then the right hand side diagram of the potential equivalence gives $\phi (1_p\cdot k) = i + \nu(k)$ for some $i \in I$. Since $\kappa$ is infinite, there exists $\alpha < \kappa$ such that $e_\alpha i = 0$. So $e_\alpha \nu(k) = e_\alpha k^p$ and $e_\alpha \phi (1_p\cdot k) = e_\alpha k^p \neq e_\alpha k$, a contradiction.  

\medskip
(ii) \emph{The field $K$ is not perfect.} Then the rings $R$ and $R(\kappa,K)$ are not isomorphic.

Assume on the contrary that there exists a ring isomorphism $\psi : R(\kappa,K) \to R$. Then for each $\alpha < \kappa$, $\psi(e_\alpha)$ is a primitive idempotent in $I$. So the map $\pi : \kappa \to \kappa$ defined by $\psi(e_\alpha) = e_{\pi(\alpha)}$ is a bijection of the set $E$ which  induces a $K$-algebra automorphism $\rho$ of $P$ defined by $\rho ((k_\alpha)_{\alpha < \kappa}) = (k_{\pi^{-1}(\alpha)})_{\alpha < \kappa}$. Notice that $\rho \restriction R(\kappa,K)$ is an automorphism of $R(\kappa,K)$. Hence $\phi = \psi (\rho \restriction R(\kappa,K))$ is a ring isomorphism of $R(\kappa,K)$ on to $R$ such that $\phi \restriction E = \hbox{id}_E$. Since $I = K^{(\kappa)} = \Soc{R(\kappa,K)} = \Soc{R}$, for each $\alpha < \kappa$, $e_\alpha K = e_\alpha R(\kappa,K) = e_\alpha R$, whence also $\phi (e_\alpha K) = e_\alpha K$. It follows that the map $\phi_\alpha : K \to K$ defined for each $k \in K$ by $e_\alpha \phi_\alpha (k) = \phi (e_\alpha k)$ is an automorphism of the field $K$. Moreover, $\phi_\alpha$ restricts to an automorphism of the subfield $\varphi(K)$ of $K$.   

For each $k \in K$, we have $\phi (1_P\cdot k) = i_k + \nu(l_k)$ for uniquely determined $i_k \in I = \Soc R$ and $l_k \in K$. So for each $\alpha < \kappa$, $e_\alpha \phi_\alpha (k) = \phi (e_\alpha k) = \phi (e_\alpha 1_P\cdot k) = e_\alpha (i_k + \nu(l_k))$. If $e_\alpha i_k = 0$, then $\phi_\alpha (k) = (l_k)^p \in \nu(K)$, whence $k \in \nu(K)$. Since the field $K$ is not perfect, there exists $k \in K \setminus \nu(K)$. By the above, $e_\alpha i_k \neq 0$ for each $\alpha < \kappa$. However, since $i_k \in I = \bigoplus_{\alpha < \kappa} e_\alpha R$, the set of all $\alpha < \kappa$ such that $e_\alpha i_k \neq 0$ is finite, a contradiction.    
\end{example}       

\begin{example}\label{non-unique} Let $p$ be a prime integer and $K$ be the field of all rational functions of one variable $x$ over $\mathbb Z_p$. Then $K$ is a countable field which is not perfect (as $x$ is not the $p$th power of any $k \in K$). Both the $K$-algebra $R(\aleph_0,K)$ from Example ~\ref{exkappa} and the ring $R$ from Example ~\ref{s-semid} are semidirect products of $I = K^{(\aleph_0)}$ and $K$. However, by Example ~\ref{s-semid}(ii), $R(\aleph_0,K)$ and $R$ are not isomorphic as rings. So they serve as examples of non-isomorphic countable semiartinian commutative regular rings of Loewy length $2$ which are factor equivalent, since they both possess the same dimension sequence $\mathcal D = \{ ( \aleph_0, \{ ( 1 , K ) \mid \beta < \aleph_0 \} ), ( 1 , ( 1 , K ) ) \}$.

Notice that the isomorphisms of rings without unit $\varphi_0$ are $K$-linear both for the $K$-algebra $R(\aleph_0,K)$ and for $R$. However, the ring isomorphism $\varphi_1$ is $K$-linear for $R(\aleph_0,K)$, but not for $R$. Thus, $R(\aleph_0,K) \in \mathfrak R _K$, but $R \notin \mathfrak R _K$,
\end{example}

\section{Multiplicative bases}\label{multb}

We will now recall an important property of $K$-algebras studied in ~\cite{BGRS}, ~\cite{CMNI}, and ~\cite{Gr} that will also hold in our setting.

\begin{definition}\label{mult} Let $K$ be a field, $R$ be a $K$-algebra, and $B$ be a $K$-linear space basis of $R$. 

Each $r \in R$ is uniquely a $K$-linear combination of elements of $B$: $r = \sum_{b \in B} b k_b$ where $k_b = 0$ for almost all $b \in B$. We let $\supp{r} = \{ b \in B \mid k_b \neq 0 \}$ and $\cs{r} = \card {(\supp{r})}$. 

We define the following properties of the basis $B$:

\begin{enumerate}
\item[(i)] If $1 \leq k < \omega$, then $B$ is a \emph{$k$-bounded basis} of $R$, if $\cs{b b^\prime} \leq k$ for all $b, b^\prime \in B$. Equivalently, $\cs{r r^\prime} \leq k \cs{r} \cs{r^\prime}$ for all $r, r^\prime \in R$. A $1$-bounded basis is also called a \emph{weak multiplicative basis}. 

$B$ is a \emph{bounded basis} of $R$, if $B$ is $k$-bounded for some $1 \leq k < \omega$. 

\item[(ii)] $B$ is a \emph{multiplicative basis} of $R$, if either $b b^\prime = 0$ or $b b^\prime \in B$, for all $b, b^\prime \in B$.
\item[(iii)] $B$ is a \emph{strong multiplicative basis} of $R$, if $b b^\prime \in B$ for all $b, b^\prime \in B$.
\end{enumerate}
\end{definition}

For all $1 \leq k < l < \omega$, we have the following relations among the properties of a $K$-basis $B$ introduced above:

\medskip
strong multiplicative $\implies$ multiplicative $\implies$ weak multiplicative $\implies$ $k$-bounded  $\implies$ $l$-bounded $\implies$ bounded.
\medskip

In particular cases, also some of the reverse implications hold true. For example, if $R$ is an integral domain, then strong multiplicative = multiplicative, and if $K = \F_2$, then multiplicative = weak multiplicative.

Notice that if $B$ is a multiplicative basis of a $K$-algebra $R$, then the $K$-algebra structure of $R$ is completely determined by the multiplicative semigroup $B \cup \{ 0 \}$. 

$K$-algebras with multiplicative bases are abundant:

\begin{example}\label{algebras} Let $K$ be a field.

(i) Let $G$ be a any group. Then $G$ is a strong multiplicative basis of the group algebra $KG$. More in general, if $M$ is a monoid then $M$ is a strong multiplicative basis of the monoid algebra $KM$. For example, the polynomial $K$-algebra $K[X]$, where $X$ is any set of commuting variables, has a strong multiplicative basis consisting of all monic monomials in the variables from $X$. 

Still more in general, if $S$ is a semigroup such that the semigroup algebra $R = KS$ is unital, then $S$ is a strong multiplicative basis of the $K$-algebra $R$. Conversely, if $B$ is a strong multiplicative basis of a (unital) $K$-algebra $R$, then $R$ is the semigroup algebra $KB$. In this case $1 = \sum_{i < m} b_i k_i$ for some $1 \leq m < \omega$ and $b_i \in B$, $k_i \in K$ ($i < m$). Let $b \in B$. Then $b = \sum_{i < m} b b_i k_i$, whence $\sum_{i < m} k_i = 1$. Thus the map $\varphi : R \to K$ defined by $\varphi (\sum_{j < p} b_j k_j) = \sum_{j < p} k_j$ is a homomorphism of (unital) rings. In particular, $I_a = \ker {\varphi}$ is a two-sided ideal of $R$ (called the \emph{augmentation} ideal of $KB$), and $R/I_a \cong K$. If $\card{B} \geq 2$, then $I_a$ is a non-trivial ideal of $R$ which is both a maximal right and a maximal left ideal of $R$, and by induction on $n = \cs{r}$, it is easy to see that as a $K$-linear space, $I_a$ is generated by the elements of the form $b - b^\prime$ where $b, b^\prime \in B$.     

For example, if $1 < n < \omega$, $R = K^n$, and $\{ 1_i \mid i < n \}$ denotes the canonical (multiplicative) $K$-basis of $R$, then $R = KB$ where $B = \{ 1_0 \} \cup \{ 1_0 + 1_i \mid 0 < i < n \}$ is a strong multiplicative basis of $R$, and $C = \{ 1_i \mid 0 < i < n \}$ is the $K$-basis of $I_a$.    

(ii) For each $1 \leq n < \omega$, the full matrix algebra $R = M_n(K)$ over a field $K$ has a multiplicative basis $B$ consisting of the matrix units $e_{ij}$ ($i,j < n$). If $n \geq 2$, then $B$ is not strong. 

In fact, if $n \geq 2$, then $R$ has no strong multiplicative basis: otherwise, by part (i), $R$ is a semigroup algebra with a non-trivial augmentation ideal $I_a$, in contradiction with $R = M_n(K)$ being a simple ring.  
  
(iii) If $Q$ is any quiver, then the set of all paths in $Q$ forms a multiplicative basis $B$ of the path algebra $KQ$, called the \emph{standard} basis of $KQ$. If $Q$ has at least two vertices, then the standard basis is not strong. 

However, if $Q$ has no oriented cycles (i.e., $KQ$ is finite dimensional), then $KQ$ possesses a finite strong multiplicative $K$-basis: Let $\{ e_i \mid i < m \}$ be the set of all vertices (= paths of length $0$) of $Q$. Let $P$ be the (finite) set of all paths of length $> 0$ in $Q$. W.l.o.g.\ we can assume that $e_0$ is not the target of any path from $P$. Then $B = \{ e_0 + p \mid p \in P \} \cup \{ e_0, e_0 + e_1, \dots , \sum_{i < m} e_i \}$ is a strong multiplicative basis of $KQ$.     

In particular, by \cite[Theorem VII.1.7]{ASS}, any finite dimensional indecomposable basic hereditary algebra over an algebraically closed field has a strong multiplicative basis. Moreover, by ~\cite{BGRS}, any finite dimensional algebra of finite representation type over any algebraically closed field $K$ has a multiplicative basis.

(iv) Let $R$ be a $K$-algebra with a multiplicative basis $B$, and let $I \neq R$ be an ideal of $R$. Let $C = \{ b + I \mid b \in B \setminus I \}$. Then $C$ is the multiplicative basis of the $K$-algebra $R/I$, iff $I$ is a \emph{$2$-nomial} ideal, that is, $I$ can be generated by some elements of the form $b_1 - b_2$ and $b$ where $b_1, b_2, b \in B$, cf.\ \cite[Theorem 2.3]{Gr}.  

(v) Clearly, if the $K$-dimension of $R$ is $k < \omega$, then any $K$-basis of $R$ is $k$-bounded. We refer to ~\cite{CMNI} for various examples of (possibly non-associative and non-unital) algebras possessing $1$-bounded (= weak multiplicative) bases. 

As mentioned above, if $K = \F_2$ then the notions of a weak multiplicative and a multiplicative basis coincide. However, $B = \{ 1, i \}$ is a weak multiplicative $\mathbb R$-basis of the field $\mathbb C$, but $\mathbb C$ has no multiplicative $\mathbb R$-basis. The latter fact follows from our next lemma. 
\end{example}

\begin{lemma}\label{fields} Let $K$ and $K^\prime$ be fields such that $K^\prime$ is a subfield of $K$. Then $K$ has a multiplicative $K^\prime$-basis only in the trivial case of $K^\prime = K$.
\end{lemma}
\begin{proof}
Assume that $K^\prime \subsetneq K$ and there exists a multiplicative $K^\prime$-basis $B$ of $K$. 

First, we show that $1 \in B$. Indeed, since $B$ is a $K^\prime$-basis, $1$ is a $K^\prime$-linear combination of finitely many elements of $B$, say $1 = \sum_{n < m} b_n k^\prime_n$. Then $b_0 = \sum_{n < m} b_0 b_n k^\prime_n$, so $b_0 = b_0 b_n$ for some $n < m$, because $B$ is a multiplicative basis. Then $b_n = 1 \in B$.    

Let $b \in B$. Then $b^{-1} = \sum_{n < m} b_n k^\prime_n$ for some $m > 0$, $b_n \in B$, $0 \neq k^\prime_n \in K^\prime$ ($n < m$), whence $1 = \sum_{n < m} b b_n k^\prime_n$. Since $B$ is a multiplicative basis and $1 \in B$, necessarily $b^{-1} = b_n \in B$ for some $n < m$.

For each $b \in B$, let $\pi_b : B \to B$ denote the multiplication by $b$. By the above, $\pi_b$ is a permutation of $B$, and $G = \{ \pi_b \mid b \in B \}$ is a subgroup (of cardinality $\card (B) \geq 2$) of the symmetric group on $B$. Let $K^\prime G$ denote the group algebra of $G$ over $K^\prime$. Then the map $\varphi : K \to K^\prime G$ defined by $\sum_{i < j} b_i k^\prime_{i} \mapsto \sum_{i < j} \pi_{b_i} k^\prime_{i}$ is a ring isomorphism. However, $K$ is a field, while $K^\prime G$ is not, as it contains the (non-trivial) augmentation ideal, a contradiction.     
\end{proof}

The $K$-algebras that we will consider in Sections ~\ref{sectBerg}-\ref{eklofproblem} have strong multiplicative bases. Moreover, they are subalgebras of the $K$-algebras $K^\kappa$ for an infinite cardinal $\kappa$. However, the latter algebras do not even possess any bounded bases. In order to prove this claim, we will employ the techniques developed for the case of semigroup algebras by Okni\'nski (cf.\ \cite[Chapter XVI]{O}):  

\begin{lemma}\label{Oknin} Let $K$ be a field and $R$ be a regular right self-injective $K$-algebra possessing a bounded basis. Then $R$ is completely reducible.
\end{lemma}
\begin{proof} Assume $R$ is not completely reducible. As $R$ is regular, it must contain a countably infinite set of pairwise orthogonal non-zero idempotents $\{ f_i \mid i < \omega \}$ (cf.\ \cite[Corollary 2.16]{G}). Let $\{ A_i \mid i < \omega \}$ be a partition of $\omega$ such that each of the sets $A_i$ ($i < \omega$) is infinite. For each $i < \omega$, let $J_i = \bigoplus_{j \in A_i} f_j R$. Since $R$ is regular and right self-injective, and $\{ J_i \mid i < \omega \}$ is an independent family of right ideals of $R$, \cite[Lemma 9.7]{G} yields a set of pairwise orthogonal idempotents $\{ e_i \mid i < \omega \}$ in $R$ such that $J_i$ is an essential submodule of $e_iR$ for each $i < \omega$. In particular, $e_i f_j = f_j$ for all $i < \omega$ and $j \in A_i$. 

\medskip
Notice that for each $i < \omega$, the left ideal $Re_i$ is infinite dimensional: If not, then the $K$-subspace $L_i$ of $Re_i$ generated by $\{ f_j e_i \mid j \in A_i \}$ is finite dimensional, whence $L_i = \sum_{j \in F_i} K f_j e_i$ for a finite subset $F_i \subset A_i$. Let $m \in A_i \setminus F_i$. Then $f_m e_i = \sum_{j \in F_i} k_j f_j e_i$ for some $k_j \in K$ ($j \in F_i$). So $f_m = f^2_m = f_m e_i f_m \in \sum_{j \in F_i} f_j R$, a contradiction.          

Let $B$ be a bounded $K$-basis of $R$. Let $i < \omega$. We will prove that $\sup_{r \in Re_i} \cs{r} = \omega$. Since $Re_i$ is infinite dimensional, it is easy to see that there is a sequence $\bar{y} = \{ y_n \mid n < \omega \}$ of elements of $Re_i$ and a sequence of non-empty pairwise disjoint finite subsets of $B$, $\{ S_n \mid n < \omega \}$, such that $y_0 = e_i$, $S_0 = \supp{e_i}$, and for each $n < \omega$, $S_n \subseteq \supp{y_n} \subseteq \bigcup_{k \leq n} S_k$. We will now modify the sequence $\bar{y}$ by induction on $n$ to a sequence $\bar{x} = \{ x_n \mid n < \omega \}$ of elements of $Re_i$ such that for each $n < \omega$, $S_n \subseteq \supp{x_n} \subseteq \bigcup_{k \leq n} S_k$ and moreover $\supp{x_{n}} \cap S_k \neq \emptyset$ for each $k \leq n$, whence $\cs{x_n} > n$. 

First, we let $x_0 = y_0 \, (= e_i \neq 0)$, whence $\cs{x_0} > 0$. For the inductive step, if $\supp{y_{n+1}} \cap S_k \neq \emptyset$ for all $k \leq n$, we let $x_{n+1} = y_{n+1}$. Otherwise, we take the biggest $k_0 \leq n$ such that $\supp{y_{n+1}} \cap S_{k_0} = \emptyset$. If $\supp{y_{n+1} + y_{k_0}} \cap S_k \neq \emptyset$ for all $k < k_0$, we let $x_{n+1} = y_{n+1} + y_{k_0}$. Otherwise, we take the biggest $k_1 < k_0$ such that $\supp{y_{n+1} + y_{k_0}} \cap S_{k_1} = \emptyset$. Proceeding similarly, we let $x_{n+1} = y_{n+1} + y_{k_0} + y_{k_1}$ in case $\supp{y_{n+1} + y_{k_0} + y_{k_1}} \cap S_k \neq \emptyset$ for each $k < k_1$, otherwise we take the biggest $k_2 < k_1$ such that $\supp{y_{n+1} + y_{k_0} + y_{k_1}} \cap S_{k_2} = \emptyset$, etc. In at most $n$ steps, we obtain thus $x_{n+1} \in Re_i$ such that $\supp{x_{n+1}} \cap S_k \neq \emptyset$ for all $k \leq n + 1$, whence $\cs{x_{n+1}} > n + 1$. 

Let $1 \leq \ell < \infty$ be such that $B$ is $\ell$-bounded. Let $p_i = \ell i \cs{e_i}$ and $r_i = x_{p_i} \in R e_i$. Then $\cs{r_i} > p_i$.        

\medskip
Finally, we let $I = \bigoplus_{i < \omega} e_iR$ and define $f \in \Hom RIR$ by $f(e_i) = r_i$ for each $i < \omega$. As $R$ is right self-injective, there exists $r \in R$ such that $r_i = r e_i$ for all $i < \omega$. Let $n = \cs{r}$. Then $\ell n \cs{e_n} = p_n < \cs{r_n} = \cs{r e_n}$. Since $B$ is $\ell$-bounded, $\cs{r e_n} \leq \ell \cs{r} \cs{e_n} = \ell n \cs{e_n}$, a contradiction.
\end{proof}

\begin{lemma}\label{smult} Let $K$ be a field, $0 < n < \omega$, and $( R_i \mid i < n )$ be a sequence of $K$-algebras such that $R_0$ possesses a strong multiplicative basis $B_0$, containing a non-zero idempotent $e$. Moreover, assume that each of the $K$-algebras $R_i$ ($0 < i < n$) has a multiplicative basis. Then the $K$-algebra $R = \prod_{i < n} R_i$ has a strong multiplicative basis.    
\end{lemma}
\begin{proof} Let $B_i$ be a multiplicative basis of the $K$-algebra $R_i$ ($0 < i < n$). For each $i < n$, let $\nu_i : R_i \to R$ be the canonical $K$-algebra embedding. Then $B = \nu_0(B_0) \cup \bigcup_{0 < i < n} \{ \nu_0(e) + \nu_i(b) \mid b \in B_i \}$ is easily seen to be a strong multiplicative basis of $R$. 
\end{proof}
     
Now, we can prove our claim: 

\begin {theorem}\label{fproducts} Let $\kappa > 0$ be cardinal, $K$ be a field, and let $R$ be the $K$-algebra $R = K^\kappa$. If $\kappa$ is finite, then $R$ has a strong multiplicative basis, otherwise $R$ has no bounded basis.   
\end{theorem}
\begin{proof} The claim for $\kappa$ finite is a particular instance of Lemma \ref{smult} for $n = \kappa$, $R_i = K$ ($i < n$) and $e = 1$ (see also Example \ref{algebras}(i)). 

Assume that $\kappa$ is infinite. Clearly, $R$ is regular, and being a direct product of self-injective rings, it is also self-injective. However, $R$ is not completely reducible, so by Lemma \ref{Oknin}, $R$ does not possess any bounded basis.          
\end{proof}

Similarly, we obtain another consequence of Lemma \ref{Oknin}:

\begin{theorem}\label{endomorph} Let $K$ be a field, $L$ be a right $K$-linear space of dimension $\kappa > 1$, and $R$ be the endomorphism algebra of $L$. If $\kappa$ is finite, then $R$ has a multiplicative basis, but no strong multiplicative basis. If $\kappa$ is infinite, then $R$ has no bounded basis. 
\end{theorem}  
\begin{proof} The claim for $\kappa$ finite was proved in Example \ref{algebras}(ii). If $\kappa$ is infinite, then $R$ does not possess any bounded basis by Lemma \ref{Oknin}, since $R$ is regular and right self-injective, but not completely reducible, cf.\ \cite[Theorem 9.12]{G}. 
\end{proof}

In contrast with Theorem \ref{fproducts}, our Setting ~\ref{set} yields a simple inductive procedure for constructing subalgebras of infinite products of copies of $K$ which possess multiplicative bases:  
        
\begin{lemma}\label{induct} Let $K$ be a field, $\kappa$ be an infinite cardinal, and $\mathcal R = ( R_\alpha \mid \alpha < \kappa )$ be a sequence of $K$-algebras each of which possesses a multiplicative basis. Let $R = R(\kappa,K,\mathcal R)$ be the $K$-algebra defined in Setting ~\ref{set}. Then $R$ has a multiplicative basis.    

Assume moreover that $R_0$ has a strong multiplicative basis $B$ containing a nonzero idempotent. Then $R$ has a strong multiplicative basis.
\end{lemma}
\begin{proof} For each $\alpha < \kappa$, let $B_\alpha$ be a multiplicative basis of the $K$-algebra $R_\alpha$ and $\nu_\alpha : R_\alpha \to P = \prod_{\alpha < \kappa} R_\alpha$ be the canonical non-unital $K$-algebra embedding. Then $B = \bigcup_{\alpha < \kappa} \nu_\alpha(B_\alpha) \cup \{ 1_P \}$ is a multiplicative basis of $R$.    

Assume that $B_0$ is strong and there is $e \in B_0$ such that $0 \neq e = e^2$. For each $0 < \alpha < \kappa$, let $B^\prime_\alpha = \{ \nu_0(e) + \nu_\alpha(b) \mid b \in B_\alpha \}$. Then $B^\prime = \nu_0(B_0) \cup \bigcup_{0 < \alpha < \kappa} B^\prime_\alpha \cup \{ 1_P \}$ is a strong multiplicative basis of $R$. 
\end{proof}

Before proceeding, we need to recall basic properties of the $K$-algebras $R(\kappa,K,\mathcal R)$ from Setting ~\ref{set} proved in ~\cite{EGT}.

\begin{lemma}\label{forum} (\cite[Lemma 2.4]{EGT}) Let $K$ be a field and $\kappa$ an infinite cardinal. Let $\mathcal R = ( R_\alpha \mid \alpha < \kappa )$ be a sequence of $K$-algebras and $R(\kappa,K,\mathcal R)$ be the $K$-algebra defined in Setting ~\ref{set}. Then the following holds:
\begin{enumerate}
\item[(i)] If all the $R_\alpha$ ($\alpha < \kappa$) are right semiartinian of Loewy length $\sigma_\alpha$ and $\tau = \sup_{\alpha < \kappa} \sigma_\alpha$, then $R(\kappa,K,\mathcal R)$ is right semiartinian of Loewy length $\tau$ or $\tau + 1$.   
\item[(ii)] If all the $R_\alpha$ ($\alpha < \kappa$) are commutative and regular, then so is $R(\kappa,K,\mathcal R)$.
\end{enumerate}
\end{lemma}

\begin{remark}\label{comments} (i) If $\kappa$ is countable and each of the $K$-algebras $R_\alpha$ $( \alpha < \kappa )$ is regular and of countable Loewy length, and all of its layers are countably generated, then the same holds for $R(\kappa,K,\mathcal R)$. In particular, each ideal of $R(\kappa,K,\mathcal R)$ is countably generated, whence $R(\kappa,K,\mathcal R)$ is a hereditary (semiartinian regular) ring.
     
(ii) By Remark \ref{comments1}, the claims (i) and (ii) of Lemma \ref{forum} will still be true with $\alpha < \kappa$ replaced by $\alpha \in A$, where $A$ is any infinite indexing set. 
\end{remark}

Moreover, we have

\begin{lemma}\label{cases} The $K$-algebras $R(\kappa,K,\mathcal R)$ from Setting ~\ref{set} have the following property:
if $R_\alpha \in \mathfrak R _K$ for each $\alpha < \kappa$, then also $R(\kappa,K,\mathcal R) \in \mathfrak R _K$.
\end{lemma}

\begin{proof} First, we need to introduce more notation: for each $\alpha < \kappa$, denote by $\nu_\alpha$ the $\alpha$-th coordinate embedding of $R_\alpha$ into $P$, and let $( S_{\alpha \beta} \mid \beta \leq \sigma_{\alpha} )$ be the socle sequence of $R_\alpha$ and $\tau = \sup_{\alpha < \kappa} \sigma_\alpha$. For $\sigma_{\alpha} < \beta \leq \tau$, let $S_{\alpha \beta} = S_{\alpha \sigma_\alpha}$ ($= R_\alpha$). Moreover, for each $\beta \leq \tau$, let $S_\beta = \bigoplus_{\alpha < \kappa} \nu_{\alpha}(S_{\alpha \beta})$, and $S_{\tau + 1} = R(\kappa,K,\mathcal R)$.

~\cite[Lemma 2.4]{EGT} claims that the Loewy length of $R(\kappa,K,\mathcal R)$ is exactly $\tau + 1$. This is true in two cases, when

(C1) $\tau$ is not the maximum of $\sigma_\alpha$ ($\alpha < \kappa$), or

(C2) $\tau = \max_{\alpha < \kappa} \sigma_\alpha$ and the set $A = \{ \alpha < \kappa \mid \sigma_\alpha = \tau \}$ is infinite.

In both these cases, $( S_\beta \mid \beta \leq \tau + 1 )$ is the socle sequence of $R(\kappa,K,\mathcal R)$ (see ~\cite[Lemma 2.4(ii)]{EGT}). Hence, if $R_\alpha \in \mathfrak R _K$ for each $\alpha < \kappa$, then also $R(\kappa,K,\mathcal R) \in \mathfrak R _K$.

Since $S_\tau = I$ and there is a $K$-algebra isomorphism $R(\kappa,K,\mathcal R)/I \cong K$, in the cases (C1) and (C2), the dimension of the top layer $\lambda _\tau$ of $R(\kappa,K,\mathcal R)$ is $1$.

\medskip
It remains to consider the case

(C3) $\tau = \max_{\alpha < \kappa} \sigma_\alpha$ and the set $A = \{ \alpha < \kappa \mid \sigma_\alpha = \tau \}$ is finite.

Here, the Loewy length of $R(\kappa,K,\mathcal R)$ equals $\tau$. Indeed, in this case $R(\kappa,K,\mathcal R) \cong \prod_{\alpha \in A} R_{\alpha} \boxplus R(B,K,\mathcal R)$ where $B = \kappa \setminus A$ is infinite.  Moreover, since $\tau$ is not a limit ordinal, $\tau > \tau ^\prime = \sup_{\alpha \in B} \sigma_\alpha$. Replacing $\kappa$ by $B$ and $\tau$ by $\tau^\prime$, we can proceed and distinguish again three cases (C1$^\prime$), (C2$^\prime$), and (C3$^\prime$) as above. In the cases (C1$^\prime$) and (C2$^\prime$), we conclude that $R(B,K,\mathcal R)$ has Loewy length $\tau^\prime + 1$, whence $R(\kappa,K,\mathcal R)$ has Loewy length $\tau$.

In case (C3$^\prime$), we obtain a finite subset $A^\prime$ of $B$ with the infinite complement $B^\prime = B \setminus A^\prime$ and an ordinal $\tau^{\prime \prime} = \sup_{\alpha \in B^{\prime}} \sigma_\alpha < \tau ^\prime$ that yield a further decomposition $R(B,K,\mathcal R) \cong \prod_{\alpha \in A^{\prime}} R_{\alpha} \boxplus R(B^\prime,K,\mathcal R)$. Proceeding similarly, we obtain a strictly decreasing chain of ordinals $\tau > \tau^\prime > \tau^{\prime \prime} > \dots$ which has to stop. That is, the third case eventually does not occur. As the term of highest Loewy length, $\prod_{\alpha \in A} R_{\alpha}$ in the decomposition above has Loewy length $\tau$, so does $R(\kappa,K,\mathcal R)$.

It follows that in the case (C3), there is a finite subset $F \subseteq \kappa$ and a $K$-algebra decomposition $R(\kappa,K,\mathcal R) \cong \prod_{\alpha \in F} R_{\alpha} \boxplus R(G,K,\mathcal R)$ where $G = \kappa \setminus F$, the $K$-algebra $R(G,K,\mathcal R)$ is of Loewy length smaller then the Loewy length of each of the $K$-algebras $R_{\alpha}$ ($\alpha \in F$), and $R(G,K,\mathcal R)$ fits in the case (C1) or (C2). Thus, also in case (C3), if $R_\alpha \in \mathfrak R _K$ for each $\alpha < \kappa$, then also $R(\kappa,K,\mathcal R) \in \mathfrak R _K$.  
\end{proof}

\section{Conormed multiplicative bases and the algebras $B_{\alpha,n}$}\label{sectBerg}

The paper ~\cite{BGRS} actually proved a stronger result than the bare existence of a multiplicative basis $B$ for any finite dimensional algebra $R$ of finite representation type over an algebraically closed field $K$. It showed that $B$ can moreover be taken \emph{normed}, meaning that $B$ contains a complete set of pairwise orthogonal primitive idempotents of $R$ as well as a basis of each power of the Jacobson radical of $R$. 

The assumption of $K$ being algebraically closed is necessary in \cite{BGRS}: if $K$ is not algebraically closed, then there is a field $K \subseteq L$ such that $1 < [L : K ] < \aleph_0$, whence $L$ is a finite dimensional $K$-algebra of finite representation type, but $L$ has no multiplicative basis by Lemma ~\ref{fields}. 

The standard basis of the path algebra $KQ$ of any quiver $Q$ is normed. So given Gabriel's presentation of any basic indecomposable finite dimensional algebra $R$ as a factor algebra of a path algebra modulo an admissible ideal $I$ (\cite[Theorem II.3.7]{ASS}), in order to prove that $R$ has a normed multiplicative basis, it suffices to show that $I$ can be taken $2$-nomial in the sense of Example \ref{algebras}(iv).

\medskip
The rings $R$ that we consider here are regular, hence their Jacobson radical is $0$, but they are not completely reducible. Hence they contain no complete sets of pairwise orthogonal primitive idempotents. In particular, the algebras considered here have no normed bases.

However, the dual notion of a conormed basis makes perfect sense in our setting for any field $K$:

\begin{definition}\label{conormed}
Let $K$ be a field and $R$ be a semiartinian regular $K$-algebra of Loewy length $\sigma +1 \geq 2$ with the socle sequence $( S_\alpha \mid \alpha \leq \sigma + 1)$. A $K$-basis $B$ of $R$ is \emph{conormed} in case $B$ contains a $K$-basis of $S_\alpha$ for each $\alpha \leq \sigma$.  
\end{definition}

\begin{remark}\label{remcon} Let $K$ be a field and $R$ be a semiartinian regular $K$-algebra. If $R$ possesses a conormed strong multiplicative basis $B$, then $R = KB$ is a (unital) semigroup $K$-algebra, and for each $\alpha < \sigma$, there is $B_\alpha \subseteq B$ such that $S_\alpha$ is the (non-unital) semigroup $K$-algebra $KB_{\alpha}$.   
\end{remark}

First, we extend Lemma \ref{smult} to the conormed setting:

\begin{lemma}\label{smult+} Let $K$ be a field, $0 < n < \omega$, and $( R_i \mid i < n )$ be a sequence of semiartinian regular $K$-algebras such that $R_0$ possesses a conormed strong multiplicative basis, $B_0$, containing a non-zero idempotent $e \in \Soc{R_0}$. Moreover, assume that each of the $K$-algebras $R_i$ ($0 < i < n$) has a conormed multiplicative basis. Then the $K$-algebra $R = \prod_{i < n} R_i$ has a conormed strong multiplicative basis.    
\end{lemma}
\begin{proof} Let $B_i$ be a conormed multiplicative basis of the $K$-algebra $R_i$ ($0 < i < n$). For each $i < n$, let $\nu_i : R_i \to R$ be the canonical $K$-algebra embedding. Then $B = \nu_0(B_0) \cup \bigcup_{0 < i < n} \{ \nu_0(e) + \nu_i(b) \mid b \in B_i \}$ is a strong multiplicative basis of $R$ by Lemma \ref{smult}. Since $e \in \Soc{R_0}$ and $S_0 = \Soc{R} = \bigoplus_{i < n} \Soc{R_i}$, we infer that $B$ is a conormed $K$-basis of $R$.   
\end{proof}

Lemmas ~\ref{forum} and ~\ref{cases} can be employed in a recursive construction of algebras with conormed (strong) multiplicative bases.

\begin{lemma}\label{vNinduct} Let $K$ be a field and $\kappa$ an infinite cardinal. Let $\mathcal R = ( R_\beta \mid \beta < \kappa )$ be a sequence of semiartinian $K$-algebras such that for each $\beta < \kappa$, the $K$-algebra $R_\beta$ has a conormed multiplicative basis. Then the $K$-algebra $R = R(\kappa,K,\mathcal R)$ has a conormed multiplicative basis.      

Moreover, assume that the $K$-algebra $R_0$ possesses a conormed strong multiplicative basis, $B_0$, containing a non-zero idempotent $e \in \Soc{R_0}$. Then the $K$-algebra $R = R(\kappa,K,\mathcal R)$ has a conormed strong multiplicative basis.    

In particular, the $K$-algebra $R(\kappa,K)$ from Example ~\ref{exkappa} has a conormed strong multiplicative basis.
\end{lemma}

\begin{proof} We will use the notation of the proof of Lemma ~\ref{cases}. In particular, (C1), (C2) and (C3) will denote the three cases discussed there.

In cases (C1) and (C2), $( S_\beta \mid \beta \leq \tau + 1 )$ is the socle sequence of the $K$-algebra $R = R(\kappa,K,\mathcal R)$, and the top layer dimension of $R$ is $1$. So we can proceed as in the proof of Lemma ~\ref{induct}: For each $\alpha < \kappa$, we let $B_\alpha$ be a conormed multiplicative basis of the $K$-algebra $R_\alpha$ and $\nu_\alpha : R_\alpha \to P = \prod_{\alpha < \kappa} R_\alpha$ be the canonical non-unital $K$-algebra embedding. Then $B = \bigcup_{\alpha < \kappa} \nu_\alpha(B_\alpha) \cup \{ 1_R \}$ is a conormed multiplicative basis of $R$. 

Moreover, if $R_0$ possesses a conormed strong multiplicative basis, $B_0$, containing a non-zero idempotent $e \in \Soc{R_0}$, then  $B = \nu_0(B_0) \cup \bigcup_{0 < \alpha < \kappa} \{ \nu_0(e) + \nu_\alpha(b) \mid b \in B_\alpha \} \cup \{ 1_R \}$ is a conormed strong multiplicative basis of $R$.

In the case (C3), there is a finite subset $F \subseteq \kappa$ and a $K$-algebra decomposition $R \cong \prod_{\alpha \in F} R_{\alpha} \boxplus R(G,K,\mathcal R)$ where $G = \kappa \setminus F$, the $K$-algebra $R^\prime = R(G,K,\mathcal R)$ is of Loewy length smaller than the Loewy length of each of the $K$-algebras $R_{\alpha}$ ($\alpha \in F$), and $R^\prime$ fits in case (C1) or (C2). By the above, $R^\prime$ has a conormed multiplicative basis containing $1_{R^\prime}$. Since clearly the $K$-algebra $\prod_{\alpha \in F} R_{\alpha}$ has a conormed multiplicative basis, so does $R = R(\kappa,K,\mathcal R)$.

Moreover, assume that $R_0$ possesses a conormed strong multiplicative basis, $B_0$, containing a non-zero idempotent $e \in \Soc{R_0}$. If $0 \in G$, then the $K$-algebra $R^\prime$ possesses a conormed strong multiplicative basis by the above, and so does $R$ by Lemma \ref{smult+}. If $0 \in F$, then $\prod_{\alpha \in F} R_{\alpha}$ has conormed strong multiplicative basis by Lemma \ref{smult+}, and so does $R \cong \prod_{\alpha \in F} R_{\alpha} \oplus R^\prime$.   
\end{proof}  

\medskip
Now, we will define the $K$-algebras $B_{\alpha, n}$ by induction starting from the base case of $B_{0,1} = K$, and employing Lemma ~\ref{forum} in the induction steps.

\begin{definition}\label{defbngen'} Let $K$ be a field.
\begin{enumerate}
\item[(i)] $B_{0,1} := K$.
\item[(ii)] For a non-limit ordinal $\alpha = \beta + 1$, we let $\mathcal R = ( R_m \mid m < \aleph_0 )$ be the constant sequence of $K$-algebras $R_m = B_{\beta,1}$ for each $m < \aleph_0$. We let $B_{\alpha,1} := R(\aleph_0,K,\mathcal R)$. Notice that here we are in case (C2) from the proof of Lemma ~\ref{cases} for $A = \kappa = \aleph_0$.
\item[(iii)] For each limit ordinal $\alpha$, we put $R_\beta = B_{\beta,1}$ for each $\beta < \alpha$.  Using the notation of Setting ~\ref{set}, we define $B_{\alpha,1} := R(\alpha,K,\mathcal R)$ where $\mathcal R = ( R_\beta \mid \beta < \alpha )$. Here, we are in case (C1) from the proof of Lemma ~\ref{cases}.     
 \item[(iv)] For all $1 < n < \omega$ and all ordinals $\alpha$, we let $B_{\alpha,n} := B_{\alpha,1} \boxplus \dots \boxplus B_{\alpha,1}$ (the direct product of $n$ copies of the $K$-algebra $B_{\alpha,1}$).
\end{enumerate}
\end{definition}

Notice that the $K$-algebra $B_{1,1}$ is just the $K$-algebra $R(\aleph_0,K)$ from Example ~\ref{exkappa}.   

\medskip
The next lemma describes the main structural properties of the $K$-algebras $B_{\alpha,n}$:

\begin{lemma}\label{propbn} Let $\alpha$ be an ordinal and $0 < n < \omega$. Then $B_{\alpha,n} \in \mathfrak R _K$, $B_{\alpha,n}$ has Loewy length $\alpha + 1$, and the top layer dimension $n$. Moreover, if $\alpha$ is countable, then $B_{\alpha,n}$ is of countable type.
\end{lemma}
\begin{proof} By Definition ~\ref{defbngen'}(iv), it suffices to prove the claim for $n = 1$.

The claim is trivial for $\alpha = 0$, and it holds by Example ~\ref{exkappa} for $\alpha = 1$. For $\alpha > 1$, the claim follows by induction on $\alpha$ from Definition ~\ref{defbngen'} using Lemmas ~\ref{forum} and ~\ref{cases} (the case (C2) occurs in non-limit steps, and the case (C1) in the limit ones).  

The moreover part follows from the fact that the inductive use of Definition ~\ref{defbngen'} yields, for all countable ordinals $\alpha$, $K$-algebras all of whose layers are countably generated.   
\end{proof}

\begin{lemma}\label{bnmult} Let $K$ be a field, $\alpha$ an ordinal, and $0 < n < \omega$. Then the $K$-algebra $B_{\alpha,n}$ has a conormed strong multiplicative basis.
\end{lemma}
\begin{proof} By Lemma \ref{smult+}, it suffices to prove the claim for $n = 1$. The proof is by induction on $\alpha$. If $\alpha = 1$, $B_{1,1} = R(\aleph_0,K)$, so the claim follows by Lemma ~\ref{vNinduct}. For the inductive and limit steps, we use Lemma ~\ref{vNinduct} in the setting of the cases (ii) and (iii) of Definition ~\ref{defbngen'}, respectively.
\end{proof}

\begin{remark}\label{boolean} A closer look at the proof of Lemma \ref{vNinduct} in the particular setting of the $K$-algebras $B_{\alpha,n}$ shows that the conormed strong multiplicative bases constructed for these algebras in Lemma \ref{bnmult} consist of idempotents from $K^\omega$, i.e., elements of $K^\omega$ whose all components are $0$ or $1$. In particular, if $K = \F_2$ then all the $K$-algebras $B_{\alpha,n}$ are Boolean rings.     
\end{remark}

\medskip
We turn to the uniqueness of the $K$-algebras of countable type in $\mathfrak R _K$. In order to prove the uniqueness up to a $K$-algebra isomorphisms, we will surprisingly need non-uniqueness of direct sum decompositions of the algebras $B_{\alpha,1}$:

\begin{lemma}\label{swallow}      
\begin{enumerate}
\item[(i)] Let $\alpha = \beta + 1$ be a countable non-limit ordinal. Then there is a $K$-algebra isomorphism $B_{\alpha,1} \cong B_{\beta,1} \boxplus B_{\alpha,1}$.
\item[(ii)] Let $0 \leq \gamma$ be a countable ordinal and $\alpha = \gamma + n$ for some $0 < n < \omega$. Then there are $K$-algebra isomorphisms $B_{\alpha,1} \cong B_{\gamma,1} \boxplus B_{\alpha,1}$.
\item[(iii)] Let $\alpha$ be a countable limit ordinal and $\beta < \alpha$. Then there is a $K$-algebra isomorphism $B_{\alpha,1} \cong B_{\beta,1} \boxplus B_{\alpha,1}$.
\end{enumerate}
\end{lemma}
\begin{proof} (i) By Definition ~\ref{defbngen'}(ii), $B_{\alpha,1} = R(\aleph_0,K,\mathcal R)$ where $\mathcal R$ is the constant sequence $( B_{\beta,1} \mid m < \aleph_0 )$. By Lemma ~\ref{forum}, $B_{\alpha,1}$ is the $K$-subalgebra of $\prod_{m < \omega} B_{\beta,1}$ generated by $I = \bigoplus_{m < \omega} B_{\beta,1}$ and by $1 \in \prod_{m < \omega} B_{\beta,1}$.
Since  $B_{\beta,1} \oplus I$ and $I$ are both direct sums of $\aleph_0$ copies of $B_{\beta,1}$, we can define an isomorphism of $K$-algebras without unit $\varphi : B_{\beta,1} \oplus I \to I$ that takes the $i$th copy of $B_{\beta,1}$ in $B_{\beta,1} \oplus I$ to its $i$th copy in $I$. Then $\varphi$ extends to a $K$-algebra isomorphism $\psi$ such that $\psi(1_\beta + 1_\alpha) = 1_\alpha$ where $1_\beta$ denotes the unit of the $K$-algebra $B_{\beta,1}$ and $1_\alpha$ the unit of $B_{\alpha,1}$.

(ii) follows from (i), as
$$B_{\alpha,1} \cong B_{\alpha - 1,1} \boxplus B_{\alpha,1} \cong \dots \cong B_{\gamma,1} \boxplus (B_{\gamma + 1,1} \boxplus \dots \boxplus B_{\alpha -1,1} \boxplus B_{\alpha,1}) \cong \dots \cong B_{\gamma,1} \boxplus B_{\alpha,1}.$$

(iii) By Definition ~\ref{defbngen'}(iii), $B_{\alpha,1} = R(\alpha,K,\mathcal R)$ where $\mathcal R = ( B_{\gamma,1} \mid \gamma < \alpha )$. By Corollary ~\ref{decomp}, $B_{\alpha,1} \cong R(A,K,\mathcal R) \boxplus R_{\beta}$ where $A = \kappa \setminus \{ \beta \}$. By part (i), we have $B_{\beta + 1,1} \cong B_{\beta,1} \boxplus B_{\beta + 1,1}$, whence $R(A,K,\mathcal R) \cong B_{\alpha,1}$ by
Remark ~\ref{comments1}. Hence $B_{\alpha,1} \cong B_{\beta,1} \boxplus B_{\alpha,1}$.     
\end{proof}

\begin{corollary}\label{swallow+'} Let $\beta < \alpha$ be countable ordinals. Then there is a $K$-algebra isomorphism $B_{\alpha,1} \cong B_{\beta,1} \boxplus B_{\alpha,1}$.
\end{corollary}
\begin{proof} If $\alpha$ is a limit ordinal, then the claim follows by Lemma ~\ref{swallow}(iii).

Assume $\alpha$ is non-limit, $\alpha = \gamma + n$, where $\gamma$  is a countable limit ordinal or $\gamma = 0$ and $0 < n < \omega$. If $\gamma \leq \beta < \alpha$, then $\alpha = \beta + n^\prime$ for some $0 < n^\prime < n$, and the claim follows by Lemma ~\ref{swallow}(ii). If $\beta < \gamma$, then $B_{\alpha,1} \cong B_{\gamma,1} \boxplus B_{\alpha,1}$ by Lemma ~\ref{swallow}(ii), and $B_{\gamma,1} \cong B_{\beta,1} \boxplus B_{\gamma,1}$ by Lemma ~\ref{swallow}(iii), whence $B_{\alpha,1} \cong B_{\gamma,1} \boxplus B_{\alpha,1} \cong B_{\beta,1} \boxplus (B_{\gamma,1} \boxplus B_{\alpha,1}) \cong B_{\beta,1} \boxplus B_{\alpha,1}$.
\end{proof}

\section{Structure of commutative semiartinian regular algebras of countable type}\label{sectStr}

Let $K$ be a field. In this section, we will prove that if $\sigma$ is a countable ordinal and $n > 0$ a finite number, then up to an isomorphism, $B_{\sigma, n}$ is the unique $K$-algebra of countable type whose length is $\sigma + 1$ and dimension of its top layer is $n$. Hence, factor equivalence implies $K$-algebra isomorphism for all $K$-algebras of countable type.  

\begin{theorem}\label{unique'} Let $K$ be a field, and let $R \in \mathfrak R _K$ be of countable type, of Loewy length $\sigma + 1$, and of the top layer dimension $0 < n < \omega$. Then there is a $K$-algebra isomorphism $R \cong B_{\sigma, n}$.
\end{theorem}
\begin{proof} By Lemma ~\ref{topone'}, $R = R_0 \boxplus \dots \boxplus R_{n-1}$ where for each $i < n$, $R_i \in \mathfrak R _K$ is of countable type, $R_i$ has Loewy length $\sigma + 1$ and the top layer dimension $1$. So in view of Definition ~\ref{defbngen'}(iv), it suffices to prove the claim for the case of $n = 1$. Let $\mathcal S = (S_\alpha \mid \alpha \leq \sigma + 1)$ be the socle sequence of $R$. Recall that in this case $R = S_\sigma \oplus 1_R \cdot K$ by Remark ~\ref{remdec'}. We will proceed by induction on $\sigma$.

If $\sigma = 0$, then $R \cong K$, so $R \cong B_{0,1}$.

Let $0 < \sigma$. Since $R$ is of countable type, the ideal $S_{\sigma}$ of $R$ is countably generated. By ~\cite[Propositions 2.14]{G}, there exist non-zero orthogonal idempotents $\{ f_i \mid i < \omega \}$ in $R$ such that $S_{\sigma} = \bigoplus_{i < \omega} f_i R$. Then also for each $\alpha \leq \sigma$, $S_\alpha = \bigoplus_{i < \omega} f_iS_{\alpha}$.  

Since $f_i$ is an idempotent, $f_iS_\alpha = f_iR \cap S_\alpha$, and $f_iS_{\alpha+1}/f_iS_\alpha \cong (f_iS_{\alpha+1} + S_\alpha)/S_\alpha = f_i L_\alpha$ for each $i < \omega$ and $\alpha < \sigma$. Let $\sigma_i+1$ be the least (non-limit) ordinal $\beta$ such that $f_i \in S_\beta$, that is, $(f_iS_{\beta} + S_{\beta -1})/S_{\beta - 1} \neq 0$. Then $f_iS_{\sigma_i +1} = f_iR \cap S_{\sigma_i +1} = f_iR$ and there is an isomorphism of $K$-algebras without unit $\bigoplus_{i \in A_\alpha} f_iS_{\alpha+1}/f_iS_\alpha \cong L_\alpha$, where $A_\alpha = \{ i < \omega \mid \sigma_i  \geq \alpha \}$.  

Thus for each $i < \omega$, $\{ f_iS_\alpha \mid \alpha \leq \sigma_i + 1 \}$ is the socle sequence of the $K$-algebra $f_iR$. It follows that  $f_iR \in \mathfrak R _K$ is of countable type and Loewy length $\sigma_i +1 \leq \sigma$. By the inductive premise, $f_iR \cong B_{\sigma_i ,n_i}$ for some $n_i \geq 1$.

Moreover, we can w.l.o.g. assume that $n_i=1$ for all $i < \omega$, because since $B_{\sigma_i, n_i}$ equals by definition $\boxplus _{j < n_i} B_{\sigma_i, 1}$, we can consider $n_i$ algebras $B_{\sigma_i, 1}$ instead of the algebra $B_{\sigma_i, n_i}$ and continue the proof with $S_\sigma \cong \bigoplus_{i < \omega}B_{\sigma_i, n_i} = \bigoplus _{i < \omega} \bigoplus _{j < n_i} B_{\sigma_i, 1 }$. 

Assume that $\sigma$ is a non-limit ordinal. As the layer $L_{\sigma -1} = S_{\sigma}/S_{\sigma - 1}$ is infinitely generated, the set $J = \{ j < \omega \mid f_j \in S_{\sigma} \setminus S_{\sigma - 1} \} = \{ j < \omega \mid \sigma_j + 1 = \sigma \}$ is infinite.
By Corollary ~\ref{swallow+'}, for all $k \notin J$, $B_{\sigma-1,1} \cong B_{\sigma_k,1} \boxplus B_{\sigma-1,1}$. Since the set $J$ is infinite, it follows that the ideal $S_\sigma$ of $R$ is isomorphic (as a $K$-algebra without unit) to $\bigoplus_{i < \omega} B_{\sigma -1,1}$. Thus, by Lemma ~\ref{extensions}, Remark ~\ref{comments1}, and Definition ~\ref{defbngen'}(ii), $R \cong B_{\sigma,1}$.

Assume $\sigma$ is a limit ordinal. We will show that $S_\sigma \cong \bigoplus_{\alpha < \sigma} B_{\alpha,1}$ in two steps. First we will remove possible terms whose length is repeated, and then add possible terms whose length is missing.

By Lemma ~\ref{forum} (since $\sigma_i < \sigma$ and  $\sigma$ is limit) we have $\sigma = \sup_{i < \omega} (\sigma_i+1)$. For $\alpha < \sigma$ denote by $C_\alpha$ the set $\{i < \omega \mid \sigma_i = \alpha\}$. 
Let $i  < \omega$ be such that the set $C_{ \sigma_i}$ has more than one element. Since $\sigma = \sup_{i < \omega} (\sigma_i+1)$, there is $\sigma_{i_0} > \sigma_i$ and since $\sigma_j < \sigma$ for all $j < \omega$ we can iterate this argument to achieve a 
strictly increasing infinite sequence $(\sigma_{i_j} \mid {j < \card(C_{\sigma_i})})$. Let $f$ be a bijection between $\card(C_{\sigma_i})$ and $C_{\sigma_i}$. By Corollary ~\ref{swallow+'}, for all $j < \card(C_{\sigma_i})$, $B_{\sigma_{i_j},1} \cong B_{\sigma_{f(j)},1} \boxplus B_{\sigma_{i_j},1}$.  Then $S_\sigma \cong \bigoplus_{\{\alpha < \sigma \mid \exists i\in \omega ; \alpha = \sigma_i\}} B_{\alpha,1}$ as a $K$-algebra without unit.

The set $C = \sigma \setminus \{ \sigma_i + 1 \mid i < \omega \}$ is countable. Since the set $\{ \sigma_i +1\mid i < \omega \}$ is cofinal in $\sigma$,  there is a strictly increasing function $f : C \to \omega$ such that $c < \sigma_{f(c)}+1$ for each $c \in C$. Then $B_{\sigma_{f(c)},1} \cong B_{c -1,1} \boxplus B_{\sigma_{f(c)},1}$ for each $c \in C$ by Corollary ~\ref{swallow+'}, so there is a non-unital $K$-algebra isomorphism $S_\sigma \cong \bigoplus_{\alpha < \sigma} B_{\alpha,1}$. By Lemma ~\ref{extensions}, Remark ~\ref{comments1}, and Definition ~\ref{defbngen'}(iii), we conclude that $R \cong B_{\sigma,1}$.
\end{proof}

It was shown in Example ~\ref{non-unique} that the assumption of $K$-linearity of each $\varphi_\alpha$ in definition of the class $\mathfrak{R}_K$ is necessary for existence of isomorphism of any factor equivalent rings. The following examples show us the necessity of other assumptions on $R$ in Theorem ~\ref {unique'}.

There exist factor equivalent rings which has the same dimension sequence $\mathcal D = \{ ( \lambda_\alpha, \{ ( n_{\alpha\beta} , K_{\alpha\beta} ) \mid \beta < \lambda_\alpha \} ) \mid \alpha \leq \sigma \}$ but are not isomorphic, even with $\sigma$ countable and $\lambda_\alpha = \omega$ for all $\alpha < \sigma$, if we allow more than one field in their structure.

\begin{example}
Consider the $\Q$-algebras $R = S \boxplus S \boxplus T$ and  $R^\prime= S \boxplus T \boxplus T$,  where $S=R(\aleph_0, \Q)$ and $T\subseteq R(\aleph_0, \R)$ such that $T$ consisnts of  eventually constant sequences of real numbers in which the constant is rational.
 The dimension sequence of both of these rings is $\lambda_ 0 = \aleph_0$, $K_{0,\beta}=\Q$ for infinitely many $\beta < \aleph_0$ and $K_{0,\beta}=\R$ for infinitely many  $\beta < \aleph_0$, $\lambda_1 = 3$ and $K_{1,\beta}=\Q$ $(\forall \beta < 3)$. However, they are not isomorphic, because in $R^\prime$ there is an idempotent $f$ (the unit of $T$ plus the unit of $T$), such that the dimension sequence of $f R^\prime$ is $\{ ( \aleph_0 , \{ (1,\R) \mid \beta < \aleph_0 \} ), ( 2 , \{ (1,\Q)  \} ) \}$, but in $R$  if for some idempotent $g$ the structure of $g R$ has $\lambda_1 = 2$ then $g R$ contains $\Q$ in the zero-th layer.
\end{example}

There also exist factor equivalent rings in $\mathfrak R_K$ with the same dimension sequence $\mathcal D = \{ ( \lambda_\alpha, \{ ( n_{\alpha\beta} , K_{\alpha\beta} ) \mid \beta < \lambda_\alpha \} ) \mid \alpha \leq \sigma \}$ which are not isomorphic, even with $\sigma$ countable but with some $\lambda_\alpha \neq \omega$.

\begin{example}
Let $\kappa $ be an ucountable cardinal. Then $\R$-algebras $R = R(\aleph_0, \R) \boxplus  R(\kappa, \R)$ and $R^\prime=  R(\kappa, \R) \oplus  R(\kappa, \R)$ both belong $\mathfrak R _K$ and have the same dimension sequences, namely $\{ ( \kappa , \{ (1,\R) \mid \beta < \kappa \} ), ( 2 , \{ (1,\R)  \} ) \}$. However, $R$ is not isomoprhic to $R^\prime$.

Indeed, if $\varphi : R \to R^\prime$ is an isomorphism of semiartinian regular rings, then $\varphi$ induces a bijection beween the sets of all primitive idempotents of $R$ and $R^\prime$. Hence, $\varphi$ induces by restriction an isomorphism of the socles of $R$ and $R^\prime$. By induction, it follows thus that $\varphi$ induces by restriction an isomorphism of all the corresponding terms of the socle sequences of $R$ and $R^\prime$. 

In particular, if $\varphi : R \to R^\prime$ is a $K$-algebra isomorphism, then $e^\prime = \varphi(1_0)$ is an idempotent in $R^\prime \setminus \Soc{R^\prime}$ (where $1_0$ denotes the unit of $ R(\aleph_0, \R)$) such that $e^\prime f \neq 0$ exactly for $\aleph_0$ primitive idempotents $f \in R^\prime$. However, there is no such idempotent in $R^\prime$. (Here we use the fact that the primitive idempotents in $ R(\kappa, \R)$ are exactly the sequences from $K^\kappa$ that are zero everywhere except for one term, which is $1$.)     
\end{example}

We refer to \cite{S} for an example of factor equivalent rings in $\mathfrak R_K$ with the same dimension sequence $\mathcal D = \{ ( \lambda_\alpha, \{ ( n_{\alpha\beta} , K_{\alpha\beta} ) \mid \beta < \lambda_\alpha \} ) \mid \alpha \leq \sigma \}$ with uncountable $\sigma$ that are not isomorphic, though $\lambda_\alpha = \omega$ for all $\alpha \leq \sigma$. 
 
Namely, \cite{S} proved the existence of two non-isomorphic thin-tall superatomic Boolean algebras. Boolean algebras can be viewed as Boolean rings; superatomic Boolean algebras then correspond to semiartinian Boolean rings, whence they belong to $\mathfrak R _{\F_2}$. The term {\lq}thin{\rq} refers to the parameters from the dimension sequence $\lambda_\alpha$ being $\omega$ for all $\alpha \leq \sigma$, while {\lq}tall{\rq} refers to $\sigma$ being uncountable.

In particular, there follows the existence of an $\F_2$-algebra of the form $R = S + 1_R \F_2$ such that the $\F_2$-algebra $R^\prime = S \boxplus S + 1_{S \boxplus S} \F_2$ is not isomorphic to $R$, however, $\mathcal D_R  = \mathcal D _{R^\prime} =  \{ ( \omega , \{ (1,\F_2) \mid \beta < \omega \} ) \mid \alpha < \omega_1\} \cup \{ ( 2 , \{ (1,\F_2)  \} ) \}$.

\medskip
Theorem ~\ref{unique'} and Lemma ~\ref{bnmult} have the following immediate corollary:

\begin{corollary}\label{multbs}
    Let $K$ be a field, and let $R \in \mathfrak R _K$ be of countable type, then $R$ has a conormed multiplicative basis.
\end{corollary}

However, conormed multiplicative bases exist for a number of $K$-algebras that are not of countable type. For example, all $K$-algebras $R(\kappa,K) \in \mathfrak R _K$ for $\kappa \geq \aleph_1$, as well as all $B_{\alpha,n} \in \mathfrak R _K$ for $\alpha \geq \aleph_1$ possess such bases, see Lemmas ~\ref{vNinduct}(i) and ~\ref{bnmult}(i). 

The following example shows that $\mathfrak R _K$ contains also $K$-algebras with conormed multiplicative bases that are not obtained by the construction from Setting ~\ref{set}. In fact, we do not know whether all $K$-algebras in the class $\mathfrak R _K$ have conormed bases - the question appears to be open even for $K$-algebras in $\mathfrak R _K$ of Loewy length $3$ with $\lambda_0 = \aleph_0$.    

\begin{example}\label{almostdisj} Let $K$ be a field. We will identify $\omega$ with the tree $T$ consisting of all finite sequences of natural numbers, $\mathbf{n}$, with the partial ordering defined for $\mathbf{n}, \mathbf{n^\prime} \in T$ by $\mathbf{n} \leq \mathbf{n^\prime}$, iff $\mathbf{n}$ is an initial segment of $\mathbf{n^\prime}$. For each $\mathbf{n} \in T$, $1_{\mathbf{n}} \in K^{(T)}$ will denote the characteristic function of the set of all predecessors of $\mathbf{n}$ in $T$. 

$\mathcal{A}$ will denote the family of all branches of $T$. So $\card{\mathcal{A}} = 2^\omega$, $\mathcal A = \{ A_\alpha \mid \alpha < 2^\omega \}$, and for each $\alpha < 2^\omega$, $A_\alpha = \{ \mathbf{n_i} \mid i < \omega \}$, where for each $i < \omega$, $\mathbf{n_i} \in T$ has length $i$, and $\mathbf{n_i} \leq \mathbf{n_{i+1}}$. For each $\alpha < 2^\omega$, $1_{A_\alpha} \in K^T$ will denote that characteristic function of $A_\alpha$.      

The $K$-algebra $R$ is defined as a $K$-subalgebra of $K^T$ by 
$$R = \bigoplus_{\mathbf{n} \in T} 1_{\mathbf{n}} K \oplus (\bigoplus_{\alpha < 2^\omega} 1_{A_\alpha} K) \oplus 1_{K^T} \cdot K.$$ 
Since the family $\mathcal{A}$ is almost disjoint and has cardinality $2^\omega$, we infer that $R$ is a semiartinian regular $K$-subalgebra of $K^T$ with $S_1 = \Soc{R} = \bigoplus_{\mathbf{n} \in T} 1_{\mathbf{n}} K = K^{(T)}$, $S_2 = K^{(T)} \oplus \bigoplus_{\alpha < 2^\omega} 1_{A_\alpha} K$, and $S_3 = R$. Thus $R \in \mathfrak R _K$, $R$ has Loewy length $3$, $\lambda_0 = \aleph_0$, $\lambda_1 = 2^\omega$, and $\lambda_2 = 1$. 

Finally, a conormed multiplicative basis of $R$ is given by the set $B = \{ 1_{\mathbf{n}} \mid \mathbf{n} \in T \} \cup \{ 1_{A_\alpha} \mid \alpha < 2^\omega \} \cup \{ 1_{K^T} \}$. 
\end{example}

\section{Existence of strictly $\lambda$-injective modules}\label{eklofproblem}

Modules over the $K$-algebras $R(\kappa,K)$ from Example ~\ref{exkappa} can be used to answer some long standing open questions concerning existence of $\lambda$-injective modules.

\begin{definition}\label{lambda} Let $R$ be a ring and $\lambda$ an infinite cardinal. A module $M$ is said to be \emph{$\lambda$-injective} provided that for each $< \lambda$-generated right ideal $I$ of $R$ and each $\varphi \in \Hom RIM$ there exists $\psi \in \Hom RRM$ such that $\psi \restriction I = \varphi$.
\end{definition}

For example, if $R$ is a regular ring, then each module is $\aleph_0$-injective, because each finitely generated right ideal of $R$ is a direct summand in $R$.

Following ~\cite{E}, for a ring $R$, we denote by $\gamma (R)$ be the minimal infinite cardinal $\lambda$ such that each right ideal of $R$ is $< \lambda$-generated. For example, $\gamma (R) = \aleph_0$, iff $R$ is right noetherian. Notice that the classic Baer lemma says that a module $M$ is injective, iff $M$ is $\gamma(R)$-injective.

The invariant $\gamma (R)$ encodes the level of deconstructibility of classes of modules of bounded projective dimensions as follows. If $\gamma(R) = \aleph_0$ or $\gamma(R)$ is a limit cardinal, let $\lambda = \gamma(R)$. Otherwise, $\gamma(R)$ is non-limit, $\gamma(R) = \aleph_{\gamma + 1}$, and we let $\lambda = \aleph_{\gamma}$. Then for each $n < \omega$ and each module $M$ of projective dimension $\leq n$ there is an ordinal $\sigma_M$ such that $M$ is the union of a continuous strictly increasing chain of submodules $( M_\alpha \mid \alpha \leq \sigma_M )$ with $M_0 = 0$ and, for each $\alpha < \sigma_M$, the module $M_{\alpha +1}/M_\alpha$ is a $\leq \lambda$-presented and has projective dimension $\leq n$ (see e.g., ~\cite[Lemma 8.9]{GT}).
     
Given an infinite cardinal $\lambda$, one can ask for the existence of modules $M$ that are $\lambda$-injective, but not $\lambda^+$-injective. Such modules are called \emph{strictly $\lambda$-injective} in ~\cite{M}. Of course, strictly $\lambda$-injective modules cannot exist over the rings $R$ for which $\gamma (R) \leq \lambda$.  

Eklof asked in ~\cite[Remark 2(iii)]{E} whether there exist rings possessing strictly $\aleph_{\omega}$-injective modules (see also ~\cite[Problem 3.9]{M}). The renewed interest in various degrees of injectivity stems from their close relations to the number of limit models in AECs of modules discovered recently by Mazari-Armida. We refer to his paper ~\cite{M} for more details.        

The $K$-algebras $R(\kappa,K)$ from Example ~\ref{exkappa} can be used to give a positive answer to Eklof's question, and more in general, to the question of existence of strictly $\lambda$-injective modules for each infinite cardinal $\lambda$:

\begin{theorem}\label{positive} Let $\kappa$ be an infinite cardinal and $R = R(\kappa,K)$ be the $K$-subalgebra of $K^{\kappa}$ consisting of all sequences from $K^{\kappa}$ that are constant except for finitely many terms. For each cardinal $\aleph_0 \leq \lambda \leq \kappa^+$, let $M_\lambda$ denote the submodule of $K^{\kappa}$ consisting of those $r \in K^{\kappa}$ whose support has cardinality $< \lambda$. Then there is a continuous strictly increasing chain of submodules of $K^{\kappa}$
$$\Soc{R} = M_{\aleph_0} \subsetneq M_{\aleph_1} \dots \subsetneq M_{\lambda} \subsetneq M_{\lambda^+} \subsetneq \dots \subsetneq M_{\kappa^+} = K^{\kappa}.$$
Moreover, the following holds true:
\begin{enumerate}
\item[(i)] $\gamma(R) = \kappa^+$,
\item[(ii)] for each $\aleph_0 \leq \lambda \leq \kappa$, the module $M_\lambda$ is strictly $\lambda$-injective,
\item[(iii)] $M_{\kappa^+} = K^{\kappa}$ is injective.
\end{enumerate}
\end{theorem}
\begin{proof} We start by looking at the structure of ideals of $R$. For each $\alpha < \kappa$, let $e_\alpha$ denote the primitive idempotent of $R$ whose $\alpha$th term is $1$ and all other terms are $0$. Then $E = \{ e_\alpha \mid \alpha < \kappa \}$ is the set of all primitive idempotents of $R$. Clearly, for each $\alpha < \kappa$, each ideal $I$ of $R$ containing $e_\alpha$, and each homomorphism $\varphi \in \Hom {R}{I}{K^{\kappa}}$, we have either $\varphi (e_\alpha) = 0$ or else $\varphi (e_\alpha)$  has support $\{ \alpha \}$.

Let $I$ be an ideal, and $\tau$ be an infinite cardinal such that $I$ is $\leq \tau$-generated. Assume that $I \subseteq \Soc{R} (= K^{(\kappa)})$. Then there is a subset $E_I \subseteq E$ of cardinality $\leq \tau$ such that $I \subseteq J_I = \bigoplus_{e \in E_I} eR$. Since $\Soc{R}$ is a completely reducible module, the inclusions $I \subseteq J_I \subseteq \Soc{R}$ split.  

If $I \nsubseteq \Soc{R}$, there is an idempotent $f \in I \setminus \Soc{R}$. Notice that the support of $f$ is $\kappa \setminus F$ for a finite subset $F \subseteq \kappa$. Since $\Soc{R}$ is a maximal ideal in $R$, $fR + \Soc{R} = R = fR \oplus (1-f)R$. So $I = fR \oplus ((1-f)R \cap I)$. However, the support of $1-f$ is finite (equal to $F$), so
$(1-f)R \cap I$ is a finitely generated submodule of $\Soc{R}$. Thus $I$ is a finitely generated ideal in $R$, too. Since $R$ is regular, $I$ is a direct summand in $R$.    

Now, we can prove our claims: (i) By the above, if $I$ is any infinitely generated ideal of $R$, then $I$ is a direct summand in $\Soc{R} = K^{(\kappa)}$, whence $I$ is $\leq \kappa$-generated. However, $\Soc{R}$ is not $< \kappa$ generated, so $\gamma(R) = \kappa^+$.  

(ii) Assume $\aleph_0 \leq \lambda \leq \kappa^+$. Let $I$ be a $< \lambda$-generated ideal of $R$ and $\varphi \in \Hom {R}I{M_\lambda}$. If $I \nsubseteq \Soc{R}$, then $I$ is a direct summand in $R$, whence $\varphi$ extends to $R$. If $I \subseteq \Soc{R}$, then $\varphi$ extends to a $\varphi^\prime \in \Hom {R}{J_I}{M_\lambda}$. Since $J_I = \bigoplus_{e \in E_I} eR$, there exists $m \in M_\lambda$ such that $m\cdot e = \varphi^\prime (e)$ for all $e \in E_I$. Then $\varphi^\prime$ extends to $\psi \in \Hom {R}{R}{M_\kappa}$ defined by $\psi(1) = m$. Thus $M_\lambda$ is $\lambda$-injective. (Notice that the particular case of $\lambda = \kappa^+$ already gives our claim (iii), since $M_{\kappa^+}  = K^{\kappa}$.)

It remains to prove that $M_\lambda$ is not $\lambda^+$-injective. Let $I = \bigoplus_{\alpha \leq \lambda} e_\alpha R \cong K^{(\lambda)}$. Then $I$ is $\lambda$-generated, but the inclusion $I \subseteq M_{\lambda}$ does not extend to $R$, because each element of $M_{\lambda}$ has support of cardinality $< \lambda$.
\end{proof}

\medskip

\section*{Acknowledgments}

Research supported by GA\v CR 23-05148S, GAUK 101524 and SVV-2020-260721. The~authors also thank Jakub Pol\'{a}k for valuable comments.

\end{document}